\documentclass[11pt,openany,leqno]{article}
\usepackage{amsmath,amsthm,amsfonts,amssymb,amscd}
\usepackage{graphics}
\usepackage{fullpage}
\usepackage[english,frenchb]{babel}
\usepackage{enumerate}

\begin{document}

\theoremstyle{plain}
\newtheorem{thm}{\bf Th\' eor\` eme}%[section]
\newtheorem{lemm}{\bf Lemme}[thm]
\newtheorem{Propo}{\bf Proposition}[thm]
\newtheorem{coro}{\bf Corollaire}[thm]
\newtheorem{propr}{\bf Propri\'et\'e }[thm]
\newtheorem{fait}{\bf Fait}[thm]

\theoremstyle{remark}
\newtheorem{rem}{\bf Remarque$\!\!$}
\renewcommand{\therem}{}

\newcommand{\Def}{\medskip\noindent
                 {\bf Definition.}~}
\newcommand{\Exem}{\medskip\noindent
                 {\bf Example.}~}

\errorcontextlines=0

\def\sectionname{Section}
\renewcommand{\thesection}{\arabic{section}}
\renewcommand{\theequation}{\arabic{equation}}
\setcounter{section}{0}%para numerar capitulos

\def\Diff{\text{Diff}}
\def\Max{\text{max}}
\def\Log{\text{log}}
\def\loc{\text{loc}}
\def\inta{\text{int }}
\def\det{\text{det}}
\def\exp{\text{exp}}
\def\Re{\text{Re}}
\def\lip{\text{Lip}}
\def\leb{\text{Leb}}
\def\dom{\text{Dom}}
\def\diam{\text{diam}}

\def\gr{\operatorname{grau}}
\def\sen{\operatorname{sen}}
\def\ker{\operatorname{ker}}

\newcommand{\AL}{\hbox{\Large{$a$}}}
\newcommand{\AH}{\hbox{\Huge{$a$}}}
\newcommand{\loz}{{\lozenge}}

\newcommand{\Bbc}{{\mathbb{C}}}
\newcommand{\Bbf}{{\mathbb{F}}}
\newcommand{\Bbn}{{\mathbb{N}}}
\newcommand{\Bbp}{{\mathbb{P}}}
\newcommand{\Bbq}{{\mathbb{Q}}}
\newcommand{\Bbr}{{\mathbb{R}}}
\newcommand{\Bbz}{{\mathbb{Z}}}

\newcommand{\cA}{{\mathcal A}}
\newcommand{\cB}{{\mathcal{B}}}
\newcommand{\cC}{{\mathcal{C}}}
\newcommand{\cD}{{\mathcal D}}
\newcommand{\cE}{{\mathcal E}}
\newcommand{\cf}{{\mathcal F}}
\newcommand{\cG}{{\mathcal{G}}}
\newcommand{\cH}{{\mathcal H}}
\newcommand{\cI}{{\mathcal{I}}}
\newcommand{\cL}{{\mathcal{L}}}
\newcommand{\cM}{{\mathcal{M}}}
\newcommand{\cN}{{\mathcal N}}
\newcommand{\cO}{{\mathcal O}}
\newcommand{\cP}{{\mathcal{P}}}
\newcommand{\cR}{{\mathcal{R}}}
\newcommand{\cU}{{\mathcal{U}}}

\newcommand{\vep}{{\varepsilon}}
\newcommand{\al}{{\alpha}}
\newcommand{\be}{{\beta}}
\newcommand{\de}{{\delta}}
\newcommand{\De}{{\Delta}}
\newcommand{\ga}{{\gamma}}
\newcommand{\Ga}{{\Gamma}}
\newcommand{\ka}{{\kappa}}
\newcommand{\la}{{\lambda}}
\newcommand{\La}{{\Lambda}}
\newcommand{\Om}{{\Omega}}
\newcommand{\om}{{\omega}}
\newcommand{\pa}{{\partial}}
\newcommand{\ro}{{\rho}}
\newcommand{\vphi}{{\varphi}}
\newcommand{\te}{{\theta}}
\newcommand{\Te}{{\Theta}}
\newcommand{\sig}{{\sigma}}
\newcommand{\Sig}{{\Sigma}}

\newcommand{\wtA}{\widetilde{A}}
\newcommand{\wtB}{\widetilde{B}}
\newcommand{\wtcD}{\widetilde{\cD}}
\newcommand{\wtcN}{\widetilde{\cN}}
\newcommand{\wtcR}{\widetilde{\cR}}
\newcommand{\wtb}{\widetilde{b}}
\newcommand{\wtn}{\widetilde{n}}
\newcommand{\wtx}{\widetilde{x}}
\newcommand{\wtI}{\widetilde{I}}
\newcommand{\wtP}{\widetilde{P}}
\newcommand{\wtQ}{\widetilde{Q}}
\newcommand{\wtR}{\widetilde{R}}
\newcommand{\wtT}{\widetilde{T}}
\newcommand{\wtnu}{\widetilde{\nu}}
\newcommand{\wtbe}{\widetilde{\beta}}
\newcommand{\wtcdm}{\wtcD_{+}}
\newcommand{\wttm}{\wtT^{+}}
\newcommand{\rim}{R^{\infty}_{+}}
\newcommand{\wtrim}{\wtR^{\infty}_{+}}

\newcommand{\whI}{{\widehat{I}}}
\newcommand{\whd}{{\widehat{d}}}
\newcommand{\whh}{{\widehat{h}}}
\newcommand{\whm}{{\widehat{m}}}
\newcommand{\whn}{{\widehat{n}}}
\newcommand{\whs}{{\widehat{s}}}
\newcommand{\whu}{{\widehat{u}}}
\newcommand{\whz}{{\hat{z}}}
\newcommand{\whcC}{{\widehat{\cC}}}
\newcommand{\whA}{{\widehat{A}}}
\newcommand{\whB}{{\widehat{B}}}
\newcommand{\whP}{{\widehat{P}}}
\newcommand{\whR}{{\widehat{R}}}
\newcommand{\whQ}{{\widehat{Q}}}
\newcommand{\whbe}{{\widehat{\beta}}}
\newcommand{\whom}{{\widehat{\omega}}}
\newcommand{\whsig}{{\widehat{\sigma}}}

\newcommand{\ov}{\overline}
\newcommand{\ovbe}{{\overline{\beta}}}
\newcommand{\ovC}{{\overline{C}}}
\newcommand{\ovn}{{\overline{n}}}
\newcommand{\ovx}{{\overline{x}}}
\newcommand{\ovP}{{\overline{P}}}
\newcommand{\ovQ}{{\overline{Q}}}
\newcommand{\ovX}{{\overline{X}}}
\newcommand{\ovY}{{\overline{Y}}}
\newcommand{\ovfork}{{\overline{\pitchfork}}}
\newcommand{\ovforki}{{\overline{\pitchfork}_{I}}}
\newcommand{\whfork}{{\widehat{\pitchfork}}}
\newcommand{\whforki}{{\widehat{\pitchfork}_{I}}}

\newcommand{\una}{{\underline{a}}}
\newcommand{\unb}{{\underline{b}}}

\newcommand{\Lg}{\bigg\langle}
\newcommand{\Rg}{\bigg\rangle}
\newcommand{\lan}{\langle}
\newcommand{\ran}{\rangle}
\newcommand{\lgg}{{\left\langle\right.}}
\newcommand{\rg}{{\left.\right\rangle}}
\newcommand{\lV}{{\left\Vert \right.}}
\newcommand{\rV}{{\left.\right\Vert}}
\newcommand{\lv}{{\left\vert \right.}}
\newcommand{\rv}{{\left.\right\vert}}

\newcommand{{\fud}}{{\frac{1}{2}}}
\newcommand{{\fudt}}{{\tfrac{1}{2}}}
\newcommand{{\fut}}{{\frac{1}{3}}}
\newcommand{{\futt}}{{\tfrac{1}{3}}}
\newcommand{\fork}{{\pitchfork}}
\newcommand{\forki}{{\pitchfork_{I}}}
\newcommand{\cRim}{{\cR^{\infty}_{+}}}
\newcommand{\pii}{{(i)}}
\newcommand{\pjj}{{(j)}}
\newcommand{\pkk}{{(k)}}
\newcommand{\rl}{{r'}}
\newcommand{\pzz}{{(0)}}
\newcommand{\pnn}{{(n)}}
\newcommand{\prr}{{(r)}}
\newcommand{\pnu}{{(n-1)}}
\newcommand{\nuu}{{n-1}}
\newcommand{\pns}{{(n_s)}}
\newcommand{\An}{A^{(n)}}
\newcommand{\Bn}{B^{(n)}}
\newcommand{\dmu}{{\Delta^{-1}}}
\newcommand{\duz}{{\Delta^{-1}_{0}}}
\newcommand{\duu}{{\Delta^{-1}_{1}}}
\newcommand{\Ism}{{I_{m}}}
\newcommand{\Eio}{{E^{0}}}
\newcommand{\Io}{{I_{0}}}
\newcommand{\Isal}{{I_{\alpha}}}
\newcommand{\Qsal}{{Q_{\alpha}}}
\newcommand{\Psal}{{P_{\alpha}}}
\newcommand{\nsal}{{n_{\alpha}}}
\newcommand{\Isom}{{I_{\omega}}}
\newcommand{\nsom}{{n_{\omega}}}
\newcommand{\Psom}{{P_{\omega}}}
\newcommand{\Qsom}{{Q_{\omega}}}
\newcommand{\Is}{{I^{*}}}
\newcommand{\wtIs}{{\wtI^{*}}}
\newcommand{\As}{{A^{*}}}
\newcommand{\Bs}{{B^{*}}}
\newcommand{\Ps}{{P^{*}}}
\newcommand{\Qs}{{Q^{*}}}
\newcommand{\hs}{{h^{*}}}
\newcommand{\ns}{{n^{*}}}
\newcommand{\ys}{{y^{*}}}
\newcommand{\zs}{{z^{*}}}
\newcommand{\gas}{{\gamma^{*}}}
\newcommand{\kas}{{\kappa^{*}}}
\newcommand{\oms}{{\omega^{*}}}
\newcommand{\cRs}{{\cR^{*}}}
\newcommand{\nsr}{{n_{r}}}
\newcommand{\Psr}{{P_{r}}}
\newcommand{\Qsr}{{Q_{r}}}
\newcommand{\nsrl}{{n'_{r}}}
\newcommand{\Psrl}{{P'_{r}}}
\newcommand{\Qsrl}{{Q'_{r}}}
\newcommand{\nsu}{{n_{1}}}
\newcommand{\Psu}{{P_{1}}}
\newcommand{\Qsu}{{Q_{1}}}
\newcommand{\nsul}{{n'_{1}}}
\newcommand{\Psul}{{P'_{1}}}
\newcommand{\Qsul}{{Q'_{1}}}
\newcommand{\nsru}{{n_{r-1}}}
\newcommand{\Psru}{{P_{r-1}}}
\newcommand{\Qsru}{{Q_{r-1}}}
\newcommand{\nsd}{{n_{2}}}
\newcommand{\Psd}{{P_{2}}}
\newcommand{\Qsd}{{Q_{2}}}
\newcommand{\nsj}{{n_{j}}}
\newcommand{\Psj}{{P_{j}}}
\newcommand{\Qsj}{{Q_{j}}}
\newcommand{\nsm}{{n_{m}}}
\newcommand{\Psm}{{P_{m}}}
\newcommand{\Qsm}{{Q_{m}}}
\newcommand{\Esk}{{E_{k}}}
\newcommand{\nsk}{{n_{k}}}
\newcommand{\Psk}{{P_{k}}}
\newcommand{\Qsk}{{Q_{k}}}
\newcommand{\Vsk}{{V_{k}}}
\newcommand{\zk}{{z_{k}}}
\newcommand{\zo}{{z_{0}}}
\newcommand{\zell}{{z_{\ell}}}
\newcommand{\gask}{{\gamma_{k}}}
\newcommand{\vepk}{{\varepsilon_{k}}}
\newcommand{\nssk}{{n^{*}_{k}}}
\newcommand{\Pssk}{{P^{*}_{k}}}
\newcommand{\Qssk}{{Q^{*}_{k}}}
\newcommand{\Vssk}{{V^{*}_{k}}}
\newcommand{\Qsso}{{Q^{*}_{0}}}
\newcommand{\gassk}{{\gamma^{*}_{k}}}
\newcommand{\nskl}{{n'_{k}}}
\newcommand{\Pskl}{{P'_{k}}}
\newcommand{\Qskl}{{Q'_{k}}}
\newcommand{\gaskl}{{\gamma'_{k}}}
\newcommand{\Tim}{{T^{+}}}
\newcommand{\TiM}{{T^{-}}}
\newcommand{\Tsm}{{T_{+}}}
\newcommand{\TsM}{{T_{-}}}
\newcommand{\Ra}{{R_{a}}}
\newcommand{\Ral}{{R_{a'}}}
\newcommand{\Ld}{{L_{d}}}
\newcommand{\Lod}{{L^{0}_{d}}}
\newcommand{\hd}{{h_{d}}}
\newcommand{\hdl}{{h'_{d}}}
\newcommand{\Hd}{{H_{d}}}
\newcommand{\Hdl}{{H'_{d}}}
\newcommand{\Hdll}{{H''_{d}}}
\newcommand{\lad}{{\lambda_{d}}}
\newcommand{\ladl}{{\lambda'_{d}}}
\newcommand{\ladll}{{\lambda''_{d}}}
\newcommand{\mud}{{\mu_{d}}}
\newcommand{\mudl}{{\mu'_{d}}}
\newcommand{\nudl}{{\nu'_{d}}}

\newcommand{\ds}{{d_{s}}}
\newcommand{\dou}{{d^{0}_{u}}}
\newcommand{\dos}{{d^{0}_{s}}}
\newcommand{\dms}{{d^{+}_{s}}}
\newcommand{\dum}{{d^{+}_{u}}}
\newcommand{\dMs}{{d^{-}_{s}}}
\newcommand{\dss}{{d^{*}_{s}}}
\newcommand{\dsu}{{d^{*}_{u}}}

\def\za#1#2{#1^{+}_{#2}}
\def\zb#1#2{#1^{-}_{#2}}
\def\iim#1{#1^{i+1}}
\def\piim#1{#1^{(i+1)}}
\def\ism#1{#1_{i+1}}
\def\jim#1{#1^{j+1}}
\def\pism#1{#1_{(i+1)}}
\def\jiM#1{#1^{j-1}}
\def\jsM#1{#1_{j-1}}
\def\kiM#1{#1^{k-1}}
\def\ksM#1{#1_{k-1}}
\def\ksm#1{#1_{k+1}}
\def\dmuz#1{{\Delta^{-1}_{#1}}}
\def\qqim#1{{#1^{\infty}_{+}}}

%\numberwithin{section}{chapter}
%\renewcommand{\chaptername}{Chapter}
%\renewcommand{\thechapter}{\Roman{chapter}}

\title{Persistance des sous-vari\'et\'es \`a bord et \`a coins normalement dilat\'ees\\
{\large Persistence of normally expanded submanifolds with boundary or corners} \\[10pt]}
\author{Pierre BERGER}
\date{}

%%% ---------------------------------------------------------------------------------
\maketitle
\thispagestyle{empty}

\def\IMSmarkvadjust{0 pt}
\def\IMSmarkhadjust{0 pt}
\def\IMSmarkhpadding{0 pt}
\def\IMSpubltext{Published in modified form:}
\def\SBIMSMark#1#2#3{
 \font\SBF=cmss10 at 10 true pt
 \font\SBI=cmssi10 at 10 true pt
 \setbox0=\hbox{\SBF \hbox to \IMSmarkhpadding{\relax}
                Stony Brook IMS Preprint \##1}
 \setbox2=\hbox to \wd0{\hfil \SBI #2}
 \setbox4=\hbox to \wd0{\hfil \SBI #3}
 \setbox6=\hbox to \wd0{\hss
             \vbox{\hsize=\wd0 \parskip=0pt \baselineskip=10 true pt
                   \copy0 \break%
                   \copy2 \break% 
                   \copy4 \break}}
 \dimen0=\ht6   \advance\dimen0 by \vsize \advance\dimen0 by 8 true pt
                \advance\dimen0 by -\pagetotal
	        \advance\dimen0 by \IMSmarkvadjust
 \dimen2=\hsize \advance\dimen2 by .25 true in
	        \advance\dimen2 by \IMSmarkhadjust

%
%   Check for publication info
%
%  \newread\jref
  \openin2=publishd.tex
  \ifeof2\setbox0=\hbox to 0pt{}
  \else 
     \setbox0=\hbox to 3.1 true in{
                \vbox to \ht6{\hsize=3 true in \parskip=0pt  \noindent  
                {\SBI \IMSpubltext}\hfil\break
                \input publishd.tex 
                \vfill}}
  \fi
  \closein2
  \ht0=0pt \dp0=0pt
 \ht6=0pt \dp6=0pt
 \setbox8=\vbox to \dimen0{\vfill \hbox to \dimen2{\copy0 \hss \copy6}}
 \ht8=0pt \dp8=0pt \wd8=0pt
 \copy8
 \message{*** Stony Brook IMS Preprint #1, #2. #3 ***}
}

\SBIMSMark{2008/1}{March 2008}{}
%%% ----------------------------------------------------------------------

\begin{abstract}
On se propose de montrer que les vari\'et\'es \`a bord et plus g\'en\'eralement \`a coins, normalement dilat\'ees par un endomorphisme sont persistantes en tant que stratifications $a$-r\'eguli\`eres. Ce r\'esultat sera d\'emontr\'e en classe $C^s$, pour $s\ge 1$. On donne aussi un exemple simple d'une sous-vari\'et\'e \`a bord normalement dilat\'ee mais qui n'est pas persistante en tant que sous-vari\'et\'e diff\'erentiable. 
\end{abstract}
\begin{otherlanguage}{english} 
\begin{abstract}

We show that invariant submanifolds with boundary, and more generally with corners which are normally expanded by an endomorphism are persistent as $a$-regular stratifications. This result will be shown in class $C^s$, for $s\ge 1$. We present also a simple example of a submanifold with boundary which is normally expanded but non-persistent as a differentiable submanifold.
\end{abstract}
\end{otherlanguage}
\renewcommand{\theequation}{\thesection.\arabic{equation}}
\setlength{\parskip}{12pt}
\section*{Introduction}
Soit $M$ une vari\'et\'e $C^\infty$. Pour $s\ge 1$, une sous-vari\'et\'e \`a bord de $M$ de classe $C^s$ et de dimension $d$ est un sous-ensemble $N$ de $M$ tel que, pour tout point $x\in N$, il existe une carte $(U,\phi)$ d'un voisinage de $x\in M$ dans le $C^s$-atlas de $M$ engendr\'e par sa structure lisse v\'erifiant :
   \begin{equation}\label{eq1}\phi(U\cap N)=V\times \{0\},\end{equation}
o\`u $V$ est un ouvert de $\mathbb R^{d-1}\times \mathbb R^+$.

De fa\c con similaire, une sous-vari\'et\'e \`a coins de $M$ de classe $C^s$ est un sous-ensemble $N$ de $M$, tel que pour tout $x\in N$, il existe une $C^s$-carte $(U,\phi)$ d'un voisinage de $x\in M$ v\'erifiant l'\'equation (\ref{eq1}) en autorisant  $V$ \`a \^etre un ouvert de $(\mathbb R^+)^d$. 

Une {\it stratification de classe $C^s$} d'un ensemble localement compact $N$ de $M$ est une partition $\Sigma$ de $N$ en sous-vari\'et\'es de classe $C^s$ appel\'ees \emph{strates}, v\'erifiant la condition de fronti\`ere :
\[\forall (X,Y)\in \Sigma^2,\quad adh(X)\cap Y\not=\emptyset\Rightarrow adh(X)\supset Y\;\mathrm{et}\; \dim X>\dim Y.\]
Une telle stratification est \emph{$(a)$-r\'eguli\`ere} si pour toutes strates $X$ et $Y$, pour toute suite $(x_n)_n\in X^{\mathbb N}$ convergeant vers $x\in Y$, telle que $T_{x_n} X$ converge vers un certain sous-espace $P$ de $T_xM$, 
l'espace $T_xY$ est contenu dans $P$.

Une sous-vari\'et\'e \`a bord de classe $C^s$ d\'efinit canoniquement une stratification $\Sigma=(\partial N,\mathring N)$ de classe $C^s$ sur $N$ dont les strates sont le  bord $\partial N$ et l'int\'erieur $\mathring N$ de $N$. Une telle stratification est $a$-r\'eguli\`ere.

De fa\c con similaire, une sous-vari\'et\'e \`a coins de classe $C^s$ d\'efinit canoniquement une stratification ($a$-r\'eguli\`ere) $\Sigma=(X_i)_{i=1}^d$ de classe $C^s$ sur $N$ dont chaque strate $X_{k}$ est la $C^s$-vari\'et\'e form\'ee des points $x\in N$ qui ont exactement $k$ cordonn\'ees non nulles dans les cartes v\'erifiant l'\'equation (\ref{eq1}). 

Soit $f$ un {\it endomorphisme} de classe $C^s$ de $M$. Autrement dit, $f$ est une application de classe $C^s$ de $M$ dans elle-m\^eme pouvant avoir des singularit\'es et ne pas \^etre bijective. On dira que $f$ {\it pr\'eserve} une stratification $\Sigma$ d'un sous-espace $N$ de $M$, si elle envoie chaque strate de $\Sigma$ dans elle-m\^eme. On dira que la stratification $\Sigma$ est \emph{$C^s$-persistante} si toute $C^s$-perturbation $f'$ de $f$ pr\'eserve une stratification $\Sigma'$ de classe $C^s$ d'un sous-espace $N'$, telle que :
\begin{itemize} \item $N$ est hom\'eomorphe \`a $N'$ via une application $h$ $C^0$-proche de l'inclusion $N\hookrightarrow M$,
\item  la restriction de $h$ \`a chaque strate $X$ de $\Sigma$ est un diff\'eomorphisme de classe $C^s$ sur une strate de $\Sigma'$, qui est $C^s$-proche de l'inclusion canonique $X\hookrightarrow M$ pour la topologie $C^s$-compact-ouverte.
\end{itemize}
Si $\Sigma$ est $a$-r\'eguli\`ere, on dira que {\it la stratification $a$-r\'eguli\`ere $\Sigma$ est $C^s$-persistante} si la stratification $\Sigma'$ d\'efinie ci-dessus est toujours $a$-r\'eguli\`ere.

Pour $s\ge 1$, l'endomorphisme $f$ \emph{dilate $s$ fois normalement la stratification} $\Sigma$, si $f$ pr\'eserve $\Sigma$ et dilate $s$ fois normalement chacune de ses strates $X$ : \\
\indent Il existe une m\'etrique riemannienne sur $M$, un r\'eel $\lambda<1$ ainsi qu'une fonction continue et positive  $C$ sur $X$ tels que, pour tout $x\in X$, pour tous vecteurs unitaires $u\in T_xN^\bot$ et  $v\in T_xX$, on a:
\[\| p\circ T_xf^n(u)\|\ge C(x)\cdot \lambda^n \cdot \max(1,\|T_xf^n(v)\|^s),\; \forall n>0\] 
avec $p$ la projection orthogonale de $T_xM$ sur $T_xN^\bot$.

Le r\'esultat principal de cette article est le th\'eor\`eme suivant:

\begin{thm}\label{varcoin} Soient $M$ une vari\'et\'e $C^\infty$, $N$ une sous-vari\'et\'e \`a coins de classe $C^s$, pour $s>0$. Soit $f$ un endomorphisme de classe $C^s$ de $M$, pr\'eservant et dilatant $s$ fois normalement la stratification $\Sigma$ induite par $N$.

Soit $N'$ un ouvert relativement compact de $N$ dont l'adh\'erence $adh(N')$ est envoy\'ee dans $N'$ par $f$. Alors la stratification $a$-r\'eguli\`ere $\Sigma_{|N'}$ sur $N'$, dont les strates sont les intersections des strates de $\Sigma$ avec $N'$, est persistante.
\end{thm}
\begin{rem} En particulier, les sous-vari\'et\'es \`a coins compactes, de classe $C^s$, d\'efinissent une stratification $a$-r\'eguli\`ere et de classe $C^s$ qui, quand elle est dilat\'ee $s$ fois normalement, est $C^s$-persistante. \end{rem}
\begin{rem} En g\'en\'eral, les sous-vari\'et\'es \`a coins ne persistent pas en tant que sous-vari\'et\'es diff\'erentiables.\end{rem}

Nous allons d'abord rappeler l'ingr\'edient principal de la preuve de ce th\'eor\`eme : la structure de treillis de laminations et son th\'eor\`eme de persistance associ\'e. Puis, nous allons \'enoncer le th\'eor\`eme \ref{varcoin} dans le cas particulier et plus clair des vari\'et\'es \`a bord compactes de classe $C^1$. Ensuite, nous exposerons un exemple simple d'une vari\'et\'e \`a bord qui est normalement dilat\'ee mais pas persistante. Enfin, nous allons montrer ce cas particulier puis le cas g\'en\'eral du th\'eor\`eme principal de cet article. 

Ce travail n'aurait pas pu \^etre r\'ealis\'e sans la direction de J-C. Yoccoz durant ma th\`ese. 
j'exprime \'egalement ma reconnaissance P. Pansu pour de nombreuses discussions. Ce travail s'est d\'eroul\'e \`a l'universit\'e Paris Sud (Orsay) et \`a l'universit\'e d'\'etat de New York (Suny Stony Brook), je remercie ces deux institutions pour leur hospitalit\'e.
\section{Structure de treillis sur les stratifications}

Soit $\Sigma$ une $C^s$-stratification d'un sous-espace localement compact $N$ d'une vari\'et\'e $M$. Une {\it structure de treillis de laminations de classe $C^s$} sur l'espace stratifi\'e $(N,\Sigma)$ est la donn\'ee pour chaque strate $X\in \Sigma$, d'un feuilletage $\mathcal L_X$ sur un voisinage ouvert $L_X$ de $X$ dans $N$ (pour la topologie induite) qui v\'erifie les conditions suivantes :
\begin{itemize}
\item[\bf - $X$ est une feuille de $\mathcal L_X$.]
\item[\bf - $(L_X, \mathcal L_X)$ est une $C^s$-lamination :] Les feuilles de $\mathcal L_X$ sont de classe $C^s$; les petites plaques de $\mathcal L_X$ varient transversalement contin\^ument dans la topologie $C^s$.
\item[\bf - Feuilletage de laminations :] \'etant donn\'ee une strate $Y$ dont l'adh\'erence contient $X$, les petites plaques de $\mathcal L_Y$ contenues dans $L_X$ sont $C^s$-feuillet\'ees par des plaques de $\mathcal L_X$; ce feuilletage est de classe $C^s$ et varie contin\^ument transversalement aux plaques de $\mathcal L_Y$. 
\end{itemize}
Comme $\mathcal L_X$ est une lamination, la stratification $\Sigma$ est alors n\'ecessairement $a$-r\'eguli\`ere.

Nous allons d\'emontrer le th\'eor\`eme \ref{varcoin} en utilisant le th\'eor\`eme suivant, qui est une restriction du r\'esultat principal de \cite{PB1} :
\begin{thm}\label{thfonda} Soient $s\ge 1$ et  $\Sigma$ une $C^s$-stratification d'un sous-espace localement compact $N$ d'une vari\'et\'e $M$. Soit $f$ un endomorphisme de $M$ qui pr\'eserve et dilate $s$ fois normalement les strates de $\Sigma$.
Si $(N,\Sigma)$ poss\`ede une structure de treillis de classe $C^s$, satisfaisant:
\begin{enumerate}[(i)]
 \item pour chaque strate $X$ de $\Sigma$, il existe un voisinage $V_X$ de $X$ dans $N$ tel que $f$ envoie chaque plaque de $\mathcal L_X$ contenue dans $V_X$ dans une feuille de $\mathcal L_X$,
 \item Chaque feuille de $\mathcal L_X$ diff\'erente de $X$ a son image par un it\'er\'e de $f$ qui est disjoint de $V_X$.
\end{enumerate}
Soit $N'$ un ouvert relativement compact dans $N$, dont l'adh\'erence est envoy\'ee par $f$ dans lui-m\^eme. Alors, la stratification $a$-r\'eguli\`ere $\Sigma_{|N'}$ de $N'$ est $C^s$-persistante.
\end{thm}
 
\section{Vari\'et\'es \`a bord normalement dilat\'ees}
 Dans le cadre de la $C^1$-persistance des vari\'et\'es \`a bord et compactes, le th\'eor\`eme \ref{varcoin} s'\'enonce ainsi:
  \label{partvarbor}
\begin{coro}\label{varbor} Soit $N$ une sous-vari\'et\'e compacte, connexe, de classe $C^1$ et \`a bord d'une vari\'et\'e diff\'erentiable lisse $M$. Soit $f$ un endomorphisme de $M$ de classe $C^1$, pr\'eservant le bord $\partial N$ et l'int\'erieur $\mathring N$ de $N$. Si $f$ dilate (une fois) normalement $\partial N$ et $\mathring N$, alors la stratification $a$-r\'eguli\`ere $(\partial N, \mathring N)$ sur $N$ est $C^1$-persistante.\\ 
\indent Autrement dit, pour toute application $f'$ $C^1$-proche de $f$, il existe deux sous-vari\'et\'es disjointes $\partial N'$ et $\mathring N'$ telles qu'il existe un hom\'eomorphisme $h$ de $N$ sur l'union $N':=\partial N'\cup \mathring N$ v\'erifiant :\begin{itemize}
\item l'application $h$ est $C^0$-proche de l'inclusion canonique de $N$ dans $M$,
\item $f'$ envoie $\partial N'$ et $\mathring N'$ dans respectivement $\partial N'$ et $\mathring N'$,
\item la restriction de $h$ \`a $\partial N$ (resp. $\mathring N$) est un $C^1$-diff\'eomorphisme sur $\partial N'$, qui est proche de l'inclusion canonique de $\partial N$ (resp. $\mathring N$) dans $M$ pour la topologie $C^1$-compact-ouverte.
\end{itemize}
\end{coro}
\paragraph{Remarque } En g\'en\'eral $N'$ n'est pas une sous-vari\'et\'e \`a bord de classe  $C^1$ car il n'y a pas de direction transverse \`a $\partial N$ tangente \`a $N'$. La vari\'et\'e $\mathring N'$ peut s'enrouler sur $\partial N'$ et former une stratification qui n'est pas toujours $b$-r\'eguli\`ere : il existe des suites $(x_n)_n\in \mathring N'^{\mathbb N}$ et $(y_n)_n\in \partial N'^{\mathbb  N}$ telles que :
\begin{itemize}
\item $(x_n)_n$ et $(y_n)$ convergent toutes les deux vers un point $y\in \partial N'$,
\item la famille de droites\footnote{Via une carte, on peut identifier un voisinage de $x$ dans $M$ \`a un espace euclidien.}  
$((x_ny_n))_n$ converge vers une droite $D$,
\item la famille de sous-espaces $(T_{y_n}\mathring N')_n$ converge vers un certain sous-espace $P$ de $T_{y}M$,
\item  mais $D$ n'est pas contenu dans $P$.
\end{itemize} 

\subsection{Exemple d'une sous-vari\'et\'e \`a bord normalement dilat\'ee, mais non persistante en tant que sous-vari\'et\'e de classe $C^1$}  Soient $M$ le plan $\mathbb R^2$, $N$ le segment $[-1,1]\times \{0\}$ et 
\[f\; := \; (x,y)\mapsto (x^3/2+x/2,2y)\]
 qui est un diff\'eomorphisme du plan. La vari\'et\'e \`a bord $N$ est bien normalement dilat\'ee et la diff\'erentielle de $f$ sur chacune des extr\'emit\'es est une similitude de rapport $2$ .\\
\indent On perturbe maintenant $f$ au voisinage d'une des extr\'emit\'es $A$ de $N$ de fa\c con \`a ce que, sur une boule $B$ centr\'ee en $A$, la perturbation $f'$ soit \'egale \`a la composition d'une petite rotation $R$ centr\'ee en $A$ avec l'homoth\'etie $H$ centr\'ee en $A$ et de rapport 2. Le th\'eor\`eme \ref{varbor} assure l'existence d'une stratification $(N',(\partial N',\mathring N'))$ proche de $N$ et pr\'eserv\'ee par cette perturbation $f'$.\\
\indent Cette perturbation \'etant homotope \`a $f$, par une homotopie restant dans un petit voisinage de $f$ dans $C^1(M,M)$ et conservant le point fixe r\'epulsif $A$, par continuit\'e, $A$ est donc une composante connexe de $\partial N'$. On peut trouver $x\in\mathring N'\cap B$ tel que $T_x N'$ soit diff\'erent de la droite joignant $A$ \`a $x$. Sinon, au voisinage de $A$, $N'$ est une demi-droite, ce qui est absurde car la composition d'une petite rotation avec une homotopie ne pr\'eserve aucune droite. On fixe un tel $x$ et on regarde la pr\'eorbite $(x_n)_{n\le 0}$ de $R\circ H$ partant de $x$. L'application $R\circ H$ \'etant lin\'eaire et conforme, l'angle entre $T_{x_n}\mathring N$ et $\stackrel{\rightarrow}{Ax_n}$ est constant et non nul. Ainsi la stratification $(N',(\partial N',\mathring N'))$ n'est pas $b$-r\'eguli\`ere, et n'est donc pas une vari\'et\'e \`a bord.\\
\indent On remarque que ce contre-exemple peut \^etre r\'ealis\'e avec des perturbations r\'eelles analytiques de $f$ : on compose simplement $f$ par une rotation du plan centr\'ee en $A$. Pour les m\^emes raisons que ci-dessus, $A$ reste une composante connexe de $\partial N'$. On utilise alors le th\'eor\`eme de lin\'earisation de Steinberg quand l'angle de rotation est petit et irrationnel, pour conjuguer diff\'erentiablement $f'$ au voisinage de $A$ avec $R\circ H$ et ainsi se ramener au cas pr\'ec\'edent pour conclure.\\
\indent On remarque enfin que cette sous vari\'et\'e \`a bord se complexifie en une sous vari\'et\'e \`a bord de $\mathbb C^2$, dont le bord est form\'e des deux m\^emes points alors que son int\'erieur est un disque. On peut choisir cette sous-vari\'et\'e $N'$ relativement compacte ayant son adh\'erence envoy\'ee dans $N'$. On remarque que cette sous vari\'et\'e n'est persistante     qu'en tant que stratification.
\subsection {Preuve du corollaire \ref{varbor}}
Pour montrer ce corollaire, il suffit de construire une structure de treillis de laminations sur l'espace stratifi\'e $(N,(\partial N,\mathring N))$ qui v\'erifie les hypoth\`eses $(i)$ et $(ii)$ du th\'eor\`eme \ref{thfonda}.

Comme $\mathring N$ est un ouvert de $N$, on peut choisir la lamination $(L_{\mathring N},\mathcal L_{\mathring N})$ \'egale \`a la vari\'et\'e $\mathring N$ (qui est un feuilletage \`a une feuille).
 Les conditions $(i)$ et $(ii)$ associ\'ees \`a cette lamination sont alors \'evidentes.

La lamination $\mathcal L_{\partial N}$ associ\'ee \`a la strate $\partial N$ va \^etre obtenue gr\^ace \`a la construction d'une fonction r\'eelle et continue $r$ sur un voisinage $L_{\partial N}$ de $\partial N$ dans $N$ v\'erifiant les propri\'et\'es suivantes :
\begin{enumerate}
\item la pr\'eimage de 0 par $r$ est \'egale  au bord de $N$,
\item $r$ est une submersion de classe $C^1$ sur $\mathring{N}\cap L_{\partial N}$, 
\item $f$ pr\'eserve les hypersurfaces de niveau de $r$ au voisinage de $\partial N$,
\item les hypersurfaces de niveau $\lambda$ de $r$ tendent vers $\partial N$ pour la topologie $C^1$ quand $\lambda$ tend vers $0$.
\end{enumerate}
D'apr\`es 1, 2 et 4, les fibres de $r$ forment bien les feuilles d'une lamination $\mathcal L_{\partial N}$ sur $L_{\partial N}$ de classe $C^1$. D'apr\`es 1, cette lamination est bien coh\'erente avec la stratification $(\partial N,\mathring N)$. D'apr\`es 2, la condition de feuilletage de laminations est bien v\'erifi\'ee. Ainsi, le couple $((L_{\partial N},\mathcal L_{\partial N}), \mathring N)$ forme une structure de treillis sur $(\partial N,\mathring N)$. D'apr\`es 3, cette structure de treillis v\'erifie l'hypoth\`ese $(i)$ du th\'eor\`eme \ref{thfonda}. Comme $\mathcal L_{\partial N}$ est une fibration, l'hypoth\`ese $(ii)$ du th\'eor\`eme \ref{thfonda} est bien v\'erifi\'ee. Comme $f$ dilate normalement le bord et l'int\'erieur de $N$, le th\'eor\`eme \ref{thfonda} implique la $C^1$-persistance de la stratification $(\partial N,\mathring N)$.\\
\indent Pour construire la fonction $r$, on commence par mettre une structure de vari\'et\'e \`a bord $C^\infty$ sur $N$, compatible avec sa structure $C^1$ initiale\footnote{On proc\`ede comme dans  \cite{H} apr\`es avoir \'etendu $N$ en une sous-vari\'et\'e sans bord. Cependant, cette structure $C^\infty$ sera en g\'en\'eral incompatible avec la structure $C^\infty$ de $M$}. On choisit alors une m\'etrique riemannienne $g$ de classe $C^\infty$ sur $N$ et adapt\'ee\footnote{Cela signifie que la fonction $C$, dans la d\'efinition de la dilatation normale d'une sous-vari\'et\'e, peut \^etre choisie \'egale \`a 1. Une telle m\'etrique existe par la proposition 2.3 de \cite{PB1}}
 \`a la dilatation normale de $\partial N$ dans $N$. On note $\exp$ l'application exponentielle associ\'ee \`a $g$. On note $n(x)\in T_x N$ l'unique vecteur unitaire, orthogonal \`a l'espace tangent du bord de $N$ et qui pointe vers l'int\'erieur de $N$. L'application $x\mapsto n(x)$ est de classe $C^1$.\\
\indent Par compacit\'e de $\partial N$, il existe $r_0>0$ et un voisinage $V$ de $\partial N$, tels que 
\[Exp\;:\;\partial N\times [0,r_0[\rightarrow V\]
\[(x,t)\mapsto \exp_x \big(t\cdot n(x)\big)\]
soit un diff\'eomorphisme et tels que l'adh\'erence de la pr\'eimage $f^{-1}_{|N}(V)$ soit incluse dans $V$. Soient $p_1$ et $p_2$ les projections sur la premi\`ere et la deuxi\`eme coordonn\'ee de $N\times [0,r[$. On note alors $\rho$ la fonction sur $V$ \'egale \`a $p_2\circ Exp^{-1}$. C'est une submersion de $V$. Soit $\pi$ la projection de $V$ sur $\partial N$ \'egale \`a $p_1\circ Exp^{-1}$.
 Soient $t>\epsilon >0$ tels que $f^{-1}(\rho^{-1}([0,t]))$ est un compact inclus dans $\rho^{-1}([0,t-\epsilon[)$. Soit $L_{\partial N}$ l'ouvert $\rho^{-1}([0,t[)$. Soit $\phi\in C^\infty(\mathbb R)$ d\'ecroissante, valant 1 sur $]-\infty,t-\epsilon]$ et 0 sur $[t,+\infty[$. Par la suite, on s'autorisera \`a r\'eduire $t$ et donc  d'adapter $\epsilon$, $\phi$ et $L_{\partial N}$.\\
 \indent  Soient $C:= \sup_{N} \|Tf\|$ et $r'$ la fonction de classe $C^1$ sur $L_{\partial N}$ d\'efinie par :
\[	r':= (1-\phi\circ \rho)\cdot \rho+\phi\circ \rho\cdot \frac{ \rho\circ f}{C}.\]
Cela implique que le gradient $\nabla r$ de $r$ est \'egal \`a 
	\begin{equation}\label{r'}
\Rightarrow \nabla r'= (1-\phi\circ \rho)\cdot \nabla \rho+\phi\circ \rho\cdot \nabla \left(\frac{\rho\circ f}{C}\right) + \left(\frac{\rho\circ f}{C}-\rho\right)\nabla (\phi\circ \rho).\end{equation}

\indent Montrons que $r'$ est une submersion. On a $g(\nabla \rho, \nabla(\phi\circ \rho))\le 0$ et comme $C\ge \|Tf\|$, la fonction $(\rho\circ f/C-\rho)$ est n\'egative. On remarque aussi que $\nabla\rho(x)$ tend uniform\'ement vers $n\circ \pi(x)$ quand la distance entre le bord de $N$ et $x$ diminue. De plus $g(\nabla \rho,\nabla (\rho\circ f))$ est \'egal au produit scalaire de $\nabla\rho$ avec l'image par l'adjoint de $Tf$ de $\nabla\rho$. Donc par sym\'etrie de $g$, $g(\nabla \rho,\nabla (\rho\circ f))$ est \'egal \`a $g(\nabla \rho,Tf\circ \nabla \rho)$. Par cons\'equent, pour $t$ assez petit, $g(\nabla \rho,\nabla (\rho\circ f))$ est proche de $g(n\circ \pi,Tf\circ n\circ \pi)$ qui est strictement positif, par dilatation normale du bord de $N$. Ainsi, il existe $m>0$ tel que, pour $t$ assez petit et tout $x\in L_{\partial N}$, on a :
\begin{equation}\label{bordm}
g(\nabla r',\nabla \rho)> m.
\end{equation}
En particulier, $r'$ est une submersion.\\
\indent On remarque aussi que $r'=\rho$ sur $\rho^{-1}(\{t\})$ et $r'=\rho\circ f/C$ sur un voisinage de $f^{-1}\big( \rho^{-1}(\{t\})\big)$. Donc, pour $t$ assez petit, la fonction $r$ suivante est bien d\'efinie et continue :
\[r\; : \; L_{\partial N} \longrightarrow \mathbb R\]
\[x\mapsto \left\{\begin{array}{cl} 0  & \mathrm{si\;} x\in \partial N\\
 \frac{r'\circ f^n}{C^n} & \mathrm{si\;} x\in f^{-n}(L_{\partial N})\setminus f^{-n-1}(L_{\partial N}),\; n\ge 0\\
\end{array}\right.\]
\indent Une telle fonction $r$ v\'erifie donc les propri\'et\'es 1 et 3. Il ne reste donc plus qu'\`a d\'emontrer les propri\'et\'es 2 et 4 pour $t$ assez petit. On peut d\'ej\`a remarquer que la restriction  de $r$ \`a ${L_{\partial N}\setminus{\partial N}}$ est de classe $C^1$. Le reste de ces propri\'et\'es peut se d\'emontrer en utilisant des champs de c\^ones.

\indent On rappelle qu'un {\it champ de c\^ones} $\chi$ sur une partie $U$ de $N$ est un ouvert de $TN_{|U}$ tel que,  pour $x\in U$, avec l'intersection $\chi(x)$ de $\chi$ avec $T_xN$ v\'erifie :
\[\left\{\begin{array}{c}\chi(x)\not = \emptyset\\ \forall u\in \chi(x), \; \forall t\in \mathbb R\setminus \{0\},\; tu\in \chi(x)\end{array}\right.\]

 Comme le bord de $N$ est normalement dilat\'e, la propri\'et\'e 2.1.9 de \cite{PB1} entra�ne que pour tout $\epsilon>0$, il existe un champ continu de c\^ones $\chi$ sur un voisinage $U$ du bord de $\partial N$, tel que : \begin{enumerate}[(a)]
\item pour tout $x\in \partial N$, l'espace tangent $T_x\partial N$ de $\partial N$ en $x$ est maximal en tant que sous-espace vectoriel inclus dans $\chi_{x}$; tout vecteur unitaire de $\chi$ est $\epsilon$-distant d'un vecteur de $T\partial N$,
\item pour tout $x\in U$, $f^{-1}(U)$ est inclus dans $U$ et pour $x\in f^{-1}_{|N}(U)$, la pr\'eimage par $T_xf$ de $\chi_{f(x)}$ est incluse dans $\chi_x$,
\item pour $x\in f^{-1}_{|N}(U)$, l'image par $Tf$ de tout vecteur $u$ du compl\'ementaire de $\chi_x$ est un vecteur non nul (du compl\'ementaire de $\chi_{f(x)}$).\end{enumerate}
Admettons pour l'instant que toutes les lignes de niveau de $r'$ sont $C^1$-proches du bord de $N$ pour $t$ assez petit. Alors, pour $t$ assez petit, $L_{\partial N}$ est inclus dans $U$ et le noyau de $Tr'$ est inclus dans $\chi$. Par (b) et (c), la restriction de $r'\circ f^n$ \`a $f^{-1}_{|N}(L_{\partial N})$ est une submersion dont le noyau est inclus dans $\chi$. Ainsi, la restriction de $r'$ \`a $L_{\partial N}\setminus \partial N$ est une submersion dont le noyau est inclus dans $\chi$. Cela prouve la propri\'et\'e 2.\\
\indent  Pour d\'emontrer 4, on va proc\'eder par l'absurde. Soit $(x_n)_n$ une suite de points de $L_{\partial N}\setminus \partial N$ convergeant vers un certain $x\in \partial N$, telle que la suite des noyaux de $T_{x_n}r$ ne tendent pas vers $T\partial N$. Soit $P$ une valeur d'adh\'erence de cette suite. Le $d$-plan $P$ appartient donc \`a l'adh\'erence de $\chi$. Par dilatation normale, l'orbite en avant par $Tf$ d'un vecteur de $P\setminus T_x\partial N$ ne s'annule pas est tend \`a avoir un angle avec $T\partial N$ bien plus grand que celui autoris\'e par (a). Cela est contradictoire avec la propri\'et\'e 3 qui implique que l'orbite de $u$ reste dans l'adh\'erence du noyau de $Tr$ et donc dans l'adh\'erence de $\chi$.\\
\indent Il ne reste donc plus qu'\`a prouver que les noyaux de $Tr'$ peuvent \^etre choisis uniform\'ement arbitrairement proche de $T\partial N$, pour $t$ assez petit.\\
\indent On a remarqu\'e que $g(\nabla \rho,\nabla r') >0$, en (\ref{bordm}). Comme $\nabla \rho$ est proche de $n$ pour $t$ assez petit, par le th\'eor\`eme des fonctions implicites, les lignes de niveau de $r'$ sont les images par $Exp$ de sections $C^1$ du fibr\'e $\partial N\times [0,r_0[\rightarrow \partial N$. Ce fibr\'e \'etant trivial, on peut identifier de telles sections \`a des fonctions r\'eelles sur $\partial N$. Dans cette identification, la section $\sigma_\mu$ associ\'ee \`a la ligne de niveau $\mu$ de $r'$ v\'erifie :
\[r'\circ Exp(x,\sigma_\mu(x))=\mu,\; \forall x\in \partial N\]
\begin{equation}\label{2.4}\Rightarrow \partial_{T\partial N}(r'\circ Exp)(x,\sigma_\mu(x))+\partial_{\mathbb R}(r'\circ Exp)(x,\sigma_\mu(x))\cdot T\sigma_\mu(x)=0,\; \forall x\in \partial N.\end{equation}
\indent Or d'apr\`es (\ref{bordm}), on a :
\begin{equation}\label{2.4'}\partial_{\mathbb R}(r'\circ Exp)=g(\nabla r',\nabla \rho)\ge m>0.\end{equation}
Par ailleurs, on a d'apr\`es (\ref{r'}) :
\begin{equation}
T_{\partial N} (r'\circ Exp)= \frac{\phi\circ \rho}{C} \cdot T(\rho \circ f \circ Exp).
\end{equation}
La forme lin\'eaire $\partial_{T\partial N}r'$ est donc de norme inf\'erieure \`a celle de $\partial_{T\partial N}(\rho\circ f\circ Exp)$. Comme $f$ pr\'eserve le bord de $N$, la norme de $\partial_{T\partial N}(\rho\circ f\circ Exp)$ est arbitrairement petite quand $t$ tend vers 0. Par (\ref{2.4}) et (\ref{2.4'}), il en est donc de m\^eme pour  $T\sigma_\mu$.
\begin{flushright}
$\square$
\end{flushright}

\section{Vari\'et\'es \`a coins normalement dilat\'ees}
\subsection{Rappels sur les vari\'et\'es \`a coins}\label{rappel}
Pour $s\in |[1,\infty]|$, une application d'un ouvert de $\mathbb R^n_+$ dans $\mathbb R^{n'}$ est de classe $C^s$ si on peut l'\'etendre en une application de classe $C^s$ d'un ouvert de $\mathbb R^n$ dans $\mathbb R^{n'}$. La diff\'erentielle en un point d'une telle application sera la diff\'erentielle de l'une de ses extensions en ce point (qui ne d\'epend pas de l'extension). Une application d'un ouvert de $\mathbb R^n_+$ dans $\mathbb R^{n'}_+$ est de classe $C^s$ si sa composition avec l'inclusion canonique de $\mathbb R^{n'}_+$ dans $\mathbb R^{n'}$ est de classe $C^s$.\\
\indent Un $C^s$-diff\'eomorphisme d'un ouvert de $\mathbb R^n_+$ sur un ouvert de $\mathbb R^n_+$ est une application qui peut s'\'etendre en un $C^s$-diff\'eomorphisme d'un ouvert de $\mathbb  R^n$ sur un ouvert de $\mathbb  R^n$.\\ 
\indent On rappelle qu'une vari\'et\'e \`a coins $N$ de dimension $d$ est une vari\'et\'e $C^\infty$ 
model\'ee sur $\mathbb R^d_+$. Cela signifie que les changements de cartes sont des $C^\infty$-diff\'eomorphismes  
 d'ouverts de $\mathbb R^d_+$.\\
\indent Par exemple, tout cube $[0,1]^n$ pour $n\ge 0$ est muni canoniquement d'une structure de vari\'et\'e \`a coins. C'est aussi le cas du compact suivant, que l'on nome la \emph{ goutte }:
\[\big\{(x,y)\in \mathbb R^2;\; x-\sqrt{x}\le y\le \sqrt{x}-x,\; x\ge 0\big\}.\]
\indent Le {\it coindice} d'un point $x$ de $N$ est le nombre de coordonn\'ees non nulles de l'image de $x$ par une carte d'un ouvert contenant cet \'el\'ement. On note par $b N$ l'ensemble des points de $N$ d'indice sup\'erieur ou \'egal \`a $k$. 
 On note par $X_k$ l'ensemble des points de coindice $k$ ; la structure de vari\'et\'e \`a coins de $N$ induit sur $X_k$ une structure de vari\'et\'e (sans coins).
 
  Soient $x\in N$ et $E$ l'ensemble des couples $(u,\phi)$, o\`u $\phi$ est une carte de $N$ d'un voisinage de $x$ et $u$ un vecteur de $\mathbb R^n$. On d\'efinit sur $E$ une relation d'\'equivalence : deux couples $(u,\phi)$ et $(v,\psi)$ sont \'equivalents si la diff\'erentielle de $\psi\circ \phi^{-1}$ au point $\phi(x)$ envoie $u$ sur $v$. L'espace quotient est appel\'e l'{\it espace tangent} en $x$ \`a $N$. On le note $T_xN$. Par transport des structures, on obtient sur $T_xN$ une structure d'espace vectoriel de dimension $n$. 

\indent Une application continue $h$, d'une vari\'et\'e \`a coins $N$ dans une autre $N'$, est de classe $C^s$,
 si vue \`a travers des cartes $\phi$ et $\phi'$ de respectivement $N$ et $N'$, l'application $\phi'\circ h\circ  \phi^{-1}$ est de classe $C^s$ sur son ensemble de d\'efinition. Dans ce cas, pour $x\in N$, on v\'erifie que l'application $h$ induit une application lin\'eaire, dite diff\'erentielle de $h$ en $x$ et not\'ee $T_xh$, qui \`a un vecteur $v\in T_xN$ envoie la classe d'\'equivalence de  $\big(T_{\phi(x)}(\phi'\circ h\circ  \phi^{-1})(u),\phi'\big)$, o\`u $(u,\phi)$ est un repr\'esentant de $v$ et $\phi'$ une carte d'un voisinage de $h(x)$.  
L'application $h$ est une {\it immersion} (resp. {\it submersion}) si sa diff\'erentielle est injective (resp. surjective) en tout point.
 Un $C^s$-diff\'eomorphisme de vari\'et\'es \`a coins 
 est une application $C^s$ qui poss\`ede un inverse de classe $C^s$. 
 Un $C^s$-diff\'eomorphisme local (de vari\'et\'es \`a coins) est une application dont la restriction \`a un voisinage de tout point est un $C^s$-diff\'eomorphisme sur son image. 
Un {\it plongement} de classe $C^s$ est un hom\'eomorphisme sur son image, qui est aussi une immersion de classe $C^s$.

\indent On va maintenant d\'efinir une vari\'et\'e \`a coins $\partial N$ telle que $\partial N\setminus b\partial N$ s'identifie \`a $X_{d-1}$, avec $d$ la dimension de $N$. Les points de $\partial N$ sont les couples $(x,E)$ o\`u $x$ appartient \`a $bN$ et $E$ est une valeur d'adh\'erence de $(T_{x_n}X_{d-1})_n$ dans l'espace des plans de codimension 1 de $TN$, pour $(x_n)_n\in (X_{d-1})^\mathbb N$ qui tend vers $x$.\\
 \indent   Cet ensemble $\partial N$ est muni de la structure de vari\'et\'e \`a coins engendr\'ee par les cartes suivantes : pour $(x,E)\in \partial N$, on choisit une carte $\phi$ d'un voisinage distingu\'e $V$ de $x\in N$. Le sous-espace vectoriel $E$ est donc de la forme $T_{x}\phi^{-1}(\mathbb R^{k-1}\times \{0\}\times \mathbb R^{d-k})$, avec $x$ appartenant \`a $\phi^{-1}(\mathbb R^{k-1}_+\times \{0\}\times \mathbb R^{d-k}_+)$.  
On consid\`ere la restriction correspondante
\[(x,E)\mapsto\phi(x)\in \mathbb R^{k-1} \times \{0\}\times \mathbb R^{d-k}.\]
De telles applications engendrent une structure de vari\'et\'e \`a coins sur $\partial N$ de dimension $d-1$.\\
\indent La vari\'et\'e \`a coins $\partial N$ s'envoie contin\^ument dans $ N$, via l'application $p$ qui \`a $(x,E)$ associe sa premi\`ere coordonn\'ee. On remarque que $x\in X_j$ \`a exactement $(d-j)$-pr\'eimages par $p$.\\ 
\indent On appelle \emph{face de $N$} une composante connexe de $\partial N$.\\  
\begin{propr} Il existe un $C^\infty$-diff\'eomorphisme local $\phi$ de la vari\'et\'e \`a coins $\partial N\times \mathbb R^+$ sur un voisinage $V$ de $b N$ dans $N$, tel que $\phi(\cdot,0)$ est \'egal \`a $p$. L'application $\phi$ sera appel\'ee voisinage tubulaire de $\partial N$.\end{propr}
\begin{proof} La preuve d\'ecoule de la th\`ese de J. Cerf  \cite{C}. Pour toute vari\'et\'e \`a coins $N$, ce dernier construit un $C^\infty$-plongement de $N$ dans une vari\'et\'e (sans coins) de m\^eme dimension. Il construit aussi une m\'etrique riemannienne sur cette extension de $N$ telle que la vari\'et\'e $X_k$ est g\'eod\'esique, pour tout $k\ge 0$. La construction de l'application $p$ est alors classique.\end{proof}

\subsection{Preuve du r\'esultat principal (th\'eor\`eme \ref{varcoin})}
\indent Ce th\'eor\`eme se d\'emontre en construisant une structure de treillis de laminations sur $(N,\Sigma)$ v\'erifiant les conditions $(i)$ et $(ii)$ du th\'eor\`eme \ref{thfonda}. Ce dernier th\'eor\`eme implique alors la persistance de $\Sigma_{|N'}$ en tant que stratification $a$-r\'eguli\`ere. \\
\indent La construction de la structure de treillis est plus d\'elicate que celle effectu\'ee pour les vari\'et\'es \`a bord, car la dilatation normale de $X_{d-1}$ ne peut pas \^etre uniforme si $N$ n'est pas une vari\'et\'e \`a bord, avec $d$ la dimension de $N$. La m\'ethode est cependant similaire. Dans la partie \ref{I}, on va montrer qu'il suffit de construire une fonction sur $\partial N\times \mathbb R^+$ v\'erifiant des propri\'et\'es semblables \`a celles d\'ej\`a rencontr\'ees dans le cadre des vari\'et\'es \`a bord. Dans la partie \ref{II}, on construira cette fonction. On proc\`ede comme pour les vari\'et\'es \`a bord, mais par d\'efaut de dilatation normale uniforme, on est oblig\'e de changer la g\'eom\'etrie du domaine fondamental.\\

Remarquons tout d'abord que l'on peux supposer que $N$ est envoy\'ee par $f$ dans $N'$, puisque notre th\'eor\`eme concerne la persistance de la restriction de $\Sigma$ \`a $N'$ et qu'un petit voisinage de $adh(N')$ dans $N$ est envoy\'e par $f$ dans $N'$.

\indent On fixe un voisinage tubulaire $p$ de $\partial N$. On rappelle que $p$ envoie $\partial N\times \{0\}$ dans $bN$.\\
\indent Il existe $\hat V'$ et $\hat V$ deux voisinages de $\partial N\times \{0\}$ dans $\partial N\times \mathbb R^+$ et une unique application $\hat f$ de classe $C^s$ de $\hat V'$ dans $\hat V$ tel que le diagramme suivant commute : 
\[\begin{array}{rcl}\hat V'&\stackrel{\hat f}{\longrightarrow}& \hat V\\
p\downarrow&&\downarrow p\\
N&\stackrel{f}{\longrightarrow}&N\end{array}\]
On note $A_k$ l'intersection de $\partial N\times \{0\}$ avec $p^{-1}(X_k)$, pour $k\ge 1$.\\

\begin{propr}\label{4.2.2} Il existe pour chaque $k\ge 1$ un voisinage $\hat U_k$ de $A_k$ dans $\hat V$, tel que $p_{|\hat U_k}$ soit un rev\^etement \`a $d-k$ feuillets de $U_k:= p(\hat U_k)$ v\'erifiant : 
\begin{itemize}
\item $\hat f^{-1}(\hat U_k)\subset\hat U_k$ et $f^{-1}(U_k)\subset U_k$,\\

avec $F^k_x:= p_{|\hat U_k}^{-1}(x)$ pour $x\in U_k$,\\
 \item $\forall x\in f^{-1}( U_k),\quad \hat f(F_{x}^k)= F_{f(x)}^k,$
\item $\forall k\ge j,\; x\in U_k\cap U_j,\quad F_{x}^k\subset F_{x}^j$.\end{itemize}\end{propr}
La preuve de cette propri\'et\'e \'etant technique et d\'elicate, elle sera d\'emontr\'ee tout \`a la fin de la preuve du th\'eor\`eme.

\subsubsection{ \label{I} Une condition suffisante pour obtenir la persistance de la stratification} 
Pour construire une $C^s$-structure de treillis de laminations sur $\Sigma$, il suffit de trouver une fonction $r$ continue, positive, born\'ee, d\'efinie sur un voisinage ouvert $D_r$ de $\partial N\times \{0\}$ dans $\hat V'$ et v\'erifiant les propri\'et\'es suivantes :
\begin{enumerate}
\item[$P_1.$] \label{c1}il existe $C>1$ tel que  :
\[\left\{\begin{array}{cl} r\circ\hat  f= C\cdot r&\mathrm{sur}\; D_r\cap\hat  f^{-1}(D_r)\\ r^{-1}(\{0\})=\partial N\times \{0\}\end{array}\right.,\]

\item[$P_2.$] \label{c2} la restriction de $r$ \`a $D_r\setminus (\partial N\times \{0\})$ est de classe $C^s$,\\

\item[$P_3.$] \label{c3} pour $k\in \{0,\dots, d-1\}$, il existe un voisinage ouvert $L_k$ de $X_k$ dans $U_k\setminus \cup_{j<k}X_j$ tel que,
pour $x\in L_k$, la fibre $F^k_x$ est incluse dans $D_r$ et l'application  
\[r_k\;:\; x\in L_k\mapsto \big(r(y)\big)_{y\in F^k_x}\]
est localement\footnote{
Restreinte \`a un ouvert trivialisant $U$, du rev\^etement $F^k\rightarrow L_k$, l'application $r_k$ est \`a valeurs dans un espace qui s'identifie \`a $\mathbb R^{d-k}_+$. Pour cette structure, on demande que $r_{k|U}$ soit une submersion stratifi\'ee.} une submersion stratifi\'ee :
pour $l\ge k$ et $x\in X_l$, la diff\'erentielle en $x$ de $r_k$ le long de $X_l$ a un noyau de dimension minimal $k$ ; on note $T_xr_k$ cette diff\'erentielle. On demande que cette submersion stratifi\'ee v\'erifie la condition  de r\'egularit\'e $(a_f)$ de Thom : 
pour $k\le l'\le l$ et $(x_n)_n\in (X_l\cap L_k)^\mathbb N$ qui tend vers $x\in X_{l'}N\cap L_k$, le noyau de $T_{x_n}r_k$ tend vers celui de $T_xr_k$ dans la grassmanienne des $k$-plans de $TM$.
\end{enumerate}

On va montrer que l'existence d'une telle fonction est suffisante pour construire une structure de treillis qui v\'erifie les propri\'et\'es $(i)$ et $(ii)$ du th\'eor\`eme \ref{thfonda}.\\
\indent Pour ce faire, on va d\'emontrer par r\'ecurrence d\'ecroissante sur $k\in \{0,\dots , d\}$, que les fibres de $r_k$ forment une lamination $\mathcal L_k$ de classe $C^s$ qui est associ\'ee \`a la strate $X_k$ et qui v\'erifie la condition de feuilletage de laminations avec chaque lamination $(L_j,\mathcal L_j)$, pour $j>k$.\\ 
\indent Pour $k=d$, les sous ensemble $L_d$, $U_d$ et $X_d$ sont \'egaux et $F_d$ est vide. L'application $r_d$ est  nulle car a valeur dans $\mathbb R^0$. La structure de lamination $\mathcal L_d$ est donc form\'ee d'une simple feuille. On a ainsi rien \`a d\'emontrer.\\
\indent Soit $k<d$. On va montrer tout d'abord que les fibres de $r_k$ restreintes \`a $L_l\cap L_k$ forment les feuilles d'une $C^s$-lamination sur $L_l\cap L_k$ qui feuillette celle de $\mathcal L_{l|L_l\cap L_k}$, pour $l>k$.\\
\indent Par $(P_1)$, pour $j\ge l$ et $x\in L_k\cap L_l\cap X_j$, les points $r_k(x)$ et $r_l(x)$ ont chacun exactement $d-j$ coordonn\'ees nulles. Donc $(r(y))_{y\in (F_x^k\setminus F_x^l)}$ n'a aucune coordonn\'ee nulle. Comme $L_k\cap L_l$ est inclus dans $\cup_{j\ge l} X_j$, par $(P_2)$ et $(P_3)$, l'application suivante est (localement) une submersion de vari\'et\'e \`a coins de classe $C^s$ :   
\[x\in L_k\cap L_l\mapsto (r(y))_{y\in (F_x^k\setminus F_x^l)}\]
On fixe un petit voisinage distingu\'e $U$ de $x\in L_k\cap L_l$ pour la structure de vari\'et\'e \`a coins $N$. 
On identifie $U$ \`a un ouvert de $\mathbb R_+^d$ via une carte de $N$. Dans cette identification, cette submersion locale restreinte \`a $U\cap L_k\cap L_l$ peut s'\'etendre sur un ouvert de $\mathbb R^d$ et ainsi d\'efinir un feuilletage $\mathcal F$ de classe $C^s$. Par $(P_3)$, ces feuilles sont transverses \`a l'identification de la lamination $(L_l\cap L_k\cap U,\mathcal L_{l|L_k\cap L_l\cap U})$. D'apr\`es la propri\'et\'e 1.3.6 de \cite{PB1}, les intersections des feuilles de $\mathcal F$ avec celle de $\mathcal L_l$ forment les feuilles d'une $C^s$-lamination $\mathcal L_{l|U}$ qui feuillette la restriction de $\mathcal L_l$ \`a $L_k\cap L_l\cap U$.\\
\indent Comme $(L_l)_{l>k}$ est un recouvrement ouvert de $L_k\setminus X_k$, l'ensemble des fibres de $r_k$ d\'efinit les feuilles d'une $C^s$ lamination $\mathcal L_k'$ sur $L_k\setminus X_k$.\\
\indent On va montrer que l'on peut rajouter la feuille $X_k$ \`a $\mathcal L_k'$ pour former une lamination $\mathcal L_k$.\\   
\indent Pour cela, on va montrer l'existence d'un recouvrement $(U_i)_i$ de $X_k$ dans $N$ tel que l'intersection des fibres de $r_k$ avec $U_i$ sont des vari\'et\'es (des plaques) qui tendent vers $X_k\cap U_i$ dans la topologie $C^1$ compact-ouverte.\\
\indent On consid\`ere ainsi une carte $\phi\; :\; U\stackrel{\sim}{\rightarrow}\mathbb R^k\times (\mathbb R^+)^{d-k}$ d'un ouvert $U$ de $N$, intersectant $X_k$ et inclus dans $L_k$. Via $\phi$, la vari\'et\'e $U\cap X_k$ s'identifie \`a $\mathbb R^k\times \{0\}$ et la vari\'et\'e $U\cap X_l$ s'identifie \`a $\mathbb R^{k}\times (\mathbb R^+_*)^{l-k}\times \{0\}$, pour $l>k$ et avec $\mathbb R^+_*=\mathbb R^+\setminus \{0\}$.\\
\indent Comme la restriction de $r_k$ \`a $X_l$ est de classe $C^s$ et que $r_k$ y a exactement $d-l$ z\'eros, l'application suivante est bien d\'efinie et de classe $C^s$, pour $U$ assez petit et $t\in \mathbb R^{l-k}\times \{0\}$ :
\[\psi_t\; :\; (u,v)\in \mathbb R^k\times (\mathbb R^+_*)^{l-k}\mapsto (r_k\circ \phi^{-1}(u,v,0)-t)\in \mathbb R^{l-k}\times \{0\}.\]
\indent Comme $\phi$ s'annule sur $r_k^{-1}(\{t\})\cap U$ et comme, par la $(a_f)$-r\'egularit\'e de $r_k$ $(P_3)$, $\partial_v\psi_t$ est inversible, l'intersection de $r_k^{-1}(\{t\})$ avec $U$ s'identifie via $\phi$ \`a un graphe d'une fonction de $\mathbb R^k$ dans $(\mathbb R^+_*)^{l-k}\times \{0\}$. Une telle fonction est de classe $C^s$ et tend vers $0$ quand $t$ tend vers $0$, dans la topologie $C^0$ par continuit\'e de $r_k$, et dans la topologie $C^1$ par la condition $(a_f)$ v\'erifi\'ee par $r_k$ $(P_3)$.\\
\indent Comme $X_k$ est $s$ fois normalement dilat\'ee par $f$, comme $f$ pr\'eserve le feuilletage $\mathcal L_k'$ par $(P_1)$ et comme les feuilles de ce feuilletage ont un espace tangent proche de celui de $X_k$, le lemme 2.3.10 de \cite{PB1} implique que les intersections des fibres de $r_k$ avec $U$ sont des vari\'et\'es qui tendent vers $X_k\cap U$ dans la topologie $C^s$.\\
\indent De tels ouverts $U$ recouvrent $X_k$. L'union de $\mathcal L_k'$ avec $X_k$ forme donc une $C^s$-lamination $\mathcal L_k$ sur $L_k$ qui v\'erifie la condition de feuilletage. Cela ach\`eve la r\'ecurrence d\'ecroissante sur $k$.

Ainsi $(L_k,\mathcal L_k)_k$ forme une structure de treillis de laminations de classe $C^s$ sur $\Sigma$.  De plus, par la propri\'et\'e $(P_1)$, la condition $(i)$ du th\'eor\`eme est v\'erifi\'e avec $V_k:= f^{-1}(L_k)\cap L_k$, pour chaque strate $X_k$. La condition $(ii)$ du th\'eor\`eme provient elle aussi de la propri\'et\'e $(P_1)$.  
\subsubsection{\label{II} R\'ealisation de la condition suffisante }
On commence par rajouter quelques notations \`a celles d\'ej\`a \'etablies avant \ref{I}. Pour $k\ge 1$, soient $A_k:=p^{-1}(X_k)\cap \partial N\times \{0\}$ et $B_k:= adh(A_k)$. Chaque point $y\in \hat V$ poss\`ede un voisinage $U_y$ tel que $p_{|U_y}$ soit un diff\'eomorphisme sur son image. On note $p_y:=p_{|U_y}^{-1}$ l'application de $p(U_y)$ dans $U_y$.\\

\paragraph{i Construction de $R$\\}

\indent Pour $x$ appartenant \`a $\Upsilon^n:=\cap_{k=0}^n\hat  f^{{-k}}(\hat V')$, on d\'efinit :
 
\[r^n(x):= \sum_{k=0}^n p_2\circ\hat  f^k(x)\]
 o\`u $p_2$ est la projection de $\partial N\times \mathbb R^+$ sur $\mathbb R^+$.\\
 \indent  Pour $x\in \Upsilon^{n+1}$, on a
\[r^n(\hat f(x))-r^n(x)= p_2\circ\hat  f^{n+1}(x)- p_2(x).\]
Par la dilatation normale, la compacit\'e relative de $N'$ et la supposition que $N$ est envoy\'ee dans $N'$, il existe $M\ge 0$ et $T>0$ tels que, avec $R:= r^M$ restreinte \`a $\Upsilon:= R^{-1}([0,T[)$ (que l'on suppose inclus dans $\Upsilon^{M+1}$), on a :\\
\begin{itemize}
\item[i.0.] $R^{-1}(\{0\})=B_1$.\\
\item[i.1.] $\exists C>\lambda>1$ ; $\forall x\in \Upsilon$, on a $C\cdot R(x)\ge R\circ\hat  f(x)\ge \lambda\cdot R(x)$.\\
\item [i.2.] Pour tout $k\ge 0$, quitte \`a restreindre $\hat U_k$ et $U_k$, l'ouvert $\Upsilon$ contient $\hat U_k$ et l'application 

\[x\in U_k\mapsto \big(R(y)\big)_{y\in F_x^k}\] est localement une submersion de vari\'et\'es \`a coins de classe $C^s$ dans $(\mathbb R_+)^{d-k}$.\\

\end{itemize}

\paragraph{ii D\'efinition it\'erative de $r$\\}

Pour construire $r$, on va choisir un ferm\'e $U$ de $\Upsilon$, disjoint de $B_1$, tel que l'int\'erieur de $\hat f^{-1}(U)$ contient $U$ et tel que l'union $\cup_{n\ge 0} \hat f^{-n}(U)$ soit \'egale \`a $\Upsilon\setminus B_1$. On va aussi choisir une fonction $\Psi$ de classe $C^s$ sur $\Upsilon$ \'egale \`a $1$ sur $U$ et \`a $0$ sur $\Upsilon\setminus \hat f^{-1}(U)$. On pose $D:= \hat f^{-1}(U)\setminus U$ et on d\'efinit :
\[R_1:= \Psi\cdot R+(1-\Psi)\cdot \frac{R\circ \hat f}{C}\]
ainsi que :

 \[\begin{array}{rcl}
 r\; :\; \Upsilon&\longrightarrow &\qquad\mathbb R^+\\
 y&\longmapsto &\left\{
  \begin{array}{cl}
 0& \mathrm{si}\; y\in B_1\\
  \frac{R_1\circ \hat f^n(y)}{C^n}& \mathrm{si}\; y\in \hat f^{-n}(D),\; n\ge 0\\
 R_1(y) &\mathrm{si}\; y\in U
 \end{array}
 \right.\end{array}\]
Les propri\'et\'es $(P_1)$ et $(P_2)$ sont alors faciles \`a v\'erifier.\\
Pour montrer $(P_3)$, on commence par calculer le noyau de la diff\'erentielle de $r_k$ en $x\in U_k$. Pour cela, on calcule la	diff\'erentielle de $r$ en $y\in \Upsilon\setminus B_1$. On a :
\[d_yr=\left\{
\begin{array}{cl} dR_1\circ T_y\hat f^{n}& \mathrm{si}\; y\in \hat f^{-n}(D)\\
   d_yR_1& \mathrm{si}\; y\in O\end{array}\right.\]
   On note  $n_y:=0$ si $y$ appartient \`a $B_1\cup U$ et $n_y:=n$ si $y$ appartient \`a $\hat f^{-n}(D)$.\\
   \indent  On a ainsi, pour $x\in U_k$ :
   \begin{equation}\label{card1} \ker T_xr_k = \ker(d_x(R_1\circ \hat f^{n_y}\circ p_y))_{y\in F_x^k}.\end{equation}
   Mais, pour $k< d$ et $x$ appartenant \`a un petit voisinage de $X_k$, les entiers $(n_y)_{y\in F_x^k}$ n'ont aucune raison d'\^etre \'egaux. Et, comme la dilatation normale de $X_{d-1}$ n'est pas uniforme, les espaces $(\ker (d_x(R_1\circ \hat f^{n_y}\circ p_y)))_{y\in F_x^k}$ ne sont en g\'en\'eral pas proches des espaces $(\ker (d_yR))_{y\in F_x^k}$. De plus, les $n_y$ premiers it\'er\'es de $y\in F_x^k$ ne restent pas forc\'ement ni dans $\hat U_k$ ni dans un voisinage de $A_k$ o\`u sa dilatation normale agit.\\
   \indent Pour palier \`a ce probl\`eme, l'id\'ee intuitive est de regrouper par paquet les \'el\'ements de la fibre $F_x^k$, en proc\'edant par r\'ecurrence d\'ecroissante sur $k$. \\
\indent Par dilatation normale, pour chaque $k$, tout plan de dimension $k$ de $TN$, assez proche d'un plan tangent \`a $X_k$, a toutes ses pr\'eorbites par $Tf$, bas\'ees dans un certain voisinage $L_k$ de $X_k$, qui 
 tendent \`a \^etre tangentes \`a $X_k$. Appelons, de fa\c con informelle, {\it le bassin de r\'epulsion de $TX_k$} l'union de tels plans de dimension $k$ de $TN$. On va maintenant esquisser la preuve de $P_3)$, par r\'ecurrence d\'ecroissante :\\ 
\indent L'\'etape $k=d$ est toujours aussi \'evidente.\\
   \indent Pour $k<d$ et $x\in L_k$, si les entiers $(n_y)_{y\in F_x^k}$ sont tous \'egaux \`a un certain entier $n$, on s'arrange pour que, quelque soit $y\in F_x^k$, chaque point $y$ de la fibre $F^k_x$ arrive \`a $D$ en \'etant rest\'e dans $p_{|\hat U_k}^{-1}( L_k)$ et pour que $\ker( T_{f^{n} (x)} r_k)$ appartienne au bassin de r\'epulsion de $TX_k$.\\
   \indent Si les entiers $(n_y)_{y\in F_x^k}$ ne sont pas \'egaux, le minimum $n$ de cette famille n'est pas atteint pour exactement $l-k>0$ \'el\'ements de $F_x^k$. On va alors s'arranger pour que :\begin{itemize}
   \item les points $\{f^i(x)\}_{i=0}^n$ appartiennent \`a $U_k\cap L_k$,
   \item le point  $f^n(x)$ appartienne \`a $U_l$ et le noyau de $T_{f^n(x)}r_l$ intersecte le noyau de $(TR_1\circ Tp_y)_{y\in F_{f^n(x)}^k\setminus F_{f^n(x)}^l}$ en un plan de dimension $k$ qui appartient au bassin de r\'epulsion de $TX_k$.\\
   \end{itemize} 

\indent La g\'eom\'etrie de $D$ est donc dict\'ee par la dilatation normale des strates $(X_k)_k$ et par la g\'eom\'etrie des voisinages $(U_k)_k$.\\
\indent Comme la dilatation normale de ces strates n'est en g\'en\'eral uniforme que pour $k$ minimal, c'est  par r\'ecurrence croissante sur $k$ que l'on va construire $D$. On va ainsi combiner une r\'ecurrence croissante avec une r\'ecurrence d\'ecroissante...
\paragraph{iii G\'eom\'etrie du domaine fondamental \\}

On va d\'efinir dans cette partie et la suivante une famille de petits r\'eels strictement positifs $(t_k)_{k=1}^d$ par r\'ecurrence croissante : le r\'eel $t_k$ sera consid\'er\'e assez petit en fonction de $(t_j)_{j<k}$. On dira que la famille $(t_k)_{k=1}^d$ est \emph{r\'ecursivement assez petite}.\\
\indent Pour $k\in\{1,\dots ,d\}$ et $t>0$, on note :
\[W_k^t:= \left\{x\in U_k;\quad \sum_{y\in F^k_{x}} R(y)< t\right\}.\]
On pose :
\[U:=\Upsilon\setminus p_{|\hat U_k}^{-1}(W_k^{t_k})\]
Par (i.1) et la propri\'et\'e \ref{4.2.2}, pour $t<T$, on a $ f^{-1}(W_{k}^t)\subset W_k^{t/\lambda}$. On suppose donc chaque $t_k$ inf\'erieur \`a $T$, ainsi l'union $\cup_{n\ge 0} \hat f^{-n}(U)$ est \'egale \`a $\Upsilon\setminus B_1$. On suppose aussi $(t_k)_k$ r\'ecursivement assez petite, pour que $C_k:= adh(W_k\setminus \cup_{l<k} f^{-1}(W_l))$ soit un compact propre inclus dans $U_k$ et $\hat f^{-1}(\cup_{j\le k} p_{|\hat U_j}^{-1}(C_j))$ soit inclus dans l'int\'erieur de $\cup_{j\le k} p_{|\hat U_j}^{-1}(C_j)$, pour $k\in\{1,\dots, d\}$.\\
\indent On d\'emontrera \`a la fin la propri\'et\'e, non triviale, suivante :
\begin{propr}\label{proprPsi} Il existe une fonction $\Psi$ de classe $C^s$ sur $\Upsilon$, valant $1$ sur $U$ et $0$ sur $\Upsilon\setminus \hat f^{-1}(U)$ tel que, pour $(t_k)_{k=1}^d$ r\'ecursivement assez petite, le noyau $E_k(x):= \ker(dR_1\circ T_xp_y)_{y\in F_x^k}$ soit uniform\'ement proche de celui de $x\mapsto (dR\circ T_xp_y)_{y\in F_x^k}$,  pour $x\in C_k$.\end{propr}

Il s'agit maintenant de fixer $(t_k)_k$, par une r\'ecurrence croissante, en fonction de la dilatation normale. Pour convenir \`a la d\'efinition it\'erative de $r$, on va mat\'erialiser l'influence de la dilations normale des strates de $(X_k)_k$ par des champs de c\^ones.
 
\paragraph{iv Champs de c\^ones\\}

\indent La dilatation normale et la propri\'et\'e 2.1.9 de \cite{PB1} implique le 

\begin{fait}\label{faitcone} Pour $k\in\{1,\dots ,d\}$ et $\epsilon_k>0$, il existe  $t_k$ assez petit devant $(t_j)_{j<k}$ et un champ de c\^ones $\chi_k$ sur $C_k$ tels que :
\begin{enumerate}
\item  pour $x\in C_k$, $E_k(x)$ est maximal en tant que qu'espace vectoriel inclus dans $\chi_k(x)$ ; de plus, tout vecteur non nul de $\chi_k(x)$ forme un angle inf\'erieur \`a $\epsilon_k$ avec un vecteur de $E_k(x)$, 
\item le champ de c\^ones $\chi_k$ est $f_*$-stable : pour $x\in C_k\cap f^{-1}(C_k)$, le c\^one $T_xf^{-1}(\chi_k( f(x)))$ est inclus dans $\chi_k(x)$.

\end{enumerate}
\end{fait}
L'esquisse de la preuve dans ii) invite \`a fixer d\'efinitivement $(t_j)_{j=1}^d$ et $(\epsilon_j)_{j=1}^d$ tels que, pour $j\in\{1,\dots ,d-1\}$, on a :

\begin{itemize}\item[$(A_j)$] pour $i<j$ et $x\in C_j\cap C_i$, le c\^one $\chi_j(x)\cap \ker(dR_1\circ Tp_y)_{y\in F_x^i\setminus F_x^j}$ est inclus dans $\chi_i(x)$.\\\end{itemize}

Pour ce faire, on proc\`ede par r\'ecurrence croissante sur $j\in\{0,\dots ,d-1\}$ :\\
\indent L'assertion $(A_0)$ est vide de sens.\\  
\indent Soit $k\in \{1,\dots ,d-1\}$. Supposons $(\epsilon_j)_{j<k}$ et $(t_j)_{j<k}$ fix\'es pour que tout ce qui pr\'ec\`ede (et notamment les assertions $(A_j)$, pour $j<k$) soit v\'erifi\'e.\\
\indent Pour tous $i<k$ et $x\in C_k\cap C_i$, l'espace $E_i(x)$ est inclus dans le c\^one ouvert $\chi_i(x)$ et la fibre $F_x^k$ est incluse dans $F_x^i$. Donc, pour $\epsilon_k$ assez petit, le c\^one $\chi_k(x)$ qui est $\epsilon_k$-proche de $E_k(x)$, v\'erifie :
\[\chi_k(x)\cap \ker(dR_1\circ T_xp_y)_{y\in F_x^i\setminus F_x^k}\subset \chi_i(x).\]
Par compacit\'e, on peut choisir $\epsilon_k$ ind\'ependamment de $x\in C_k\cap C_i$. Ainsi, l'assertion $(A_j)$ est v\'erifi\'ee pour un tel $\epsilon_k$ que l'on fixe maintenant. On fixe aussi $t_k$ pour que tout ce qui pr\'ec\`ede soit v\'erifi\'e.

\paragraph{v V\'erification de la propri\'et\'e $(P_3)$\\}

\indent On va montrer par r\'ecurrence d\'ecroissante la propri\'et\'e suivante :   
\begin{propr}\label{terminer} Sur $C_k$, le noyau de la diff\'erentielle de $r_k$ est inclus dans $\chi_k$.\end{propr}
\begin{proof}
Pour commencer, on remarque que :
\[\emptyset =:M_{-1}\subset M_0:= p_{|\hat U_0}^{-1}(C_0)\subset \cdots M_k:= \cup_{j< k} p^{-1}_{|\hat U_j}(C_j)\subset \cdots \subset M_d=:\Upsilon\]
est une filtration. Autrement dit, la pr\'eimage par $\hat f$ de $M_k$ est incluse dans l'int\'erieur de $M_k$, pour $k\in \{0,\dots,d-1\}$.\\
\indent Comme $p^{-1}_{|\hat U_{d-1}}(C_{d-1})$ est \'egal \`a $adh(p_{|\hat U_{d-1}}^{-1}(W_{d-1})\setminus \hat f^{-1}(M_{d-2}))$ toute orbite partant de $p^{-1}_{|\hat U_{d-1}}(C_{d-1})\setminus A_{d-1}$ arrive en $D$ en \'etant rest\'ee dans $C_{d-1}$. Ainsi, le noyau de la diff\'erentielle de $r_{d-1}$ en $x\in C_{d-1}$, qui est \'egal \`a celui de $f_*^{n_y}(dR_1\circ Tp_{f^{n_y}(y)})_{y\in F^{d-1}(z)}$ par (\ref{card1}), est inclus dans $\chi_{d-1}(x)$ par les assertions 1 et 2 du fait \ref{faitcone}.\\
\indent Soit $k\in \{0,\dots,d-2\}$.  On suppose que, pour $j>k$, la propri\'et\'e \ref{terminer} est v\'erifi\'ee. Comme 
$p_{|\hat U_k}^{-1}(C_k)$ est \'egal \`a $adh(p_{|\hat U_k}^{-1}(W_k)\setminus f^{-1}(M_{k-1}))$, toute orbite partant de $p^{-1}_{|\hat U_k}(C_k)\setminus A_k$ arrive en $D$ en franchissant $(M_i)_{i\ge k}$ par ordre croissant.\\
\indent Soit $x\in C_k$. Si tous les points de $F_x^k$ arrivent en $D$ en \'etant rest\'es dans $M_k$, alors ils sont aussi tous rest\'es dans $p^{-1}_{|\hat U_k}(C_k)$. La propri\'et\'e \ref{terminer} s'obtient alors comme dans le cas $k=d-1$, car les entiers $(n_y)_{y\in F_x^k}$ sont tous \'egaux. Sinon, on consid\`ere $n\ge 0$ minimal tel qu'il existe un \'el\'ement de $y\in F_x^k$ dont l'image $\hat f^n(y)$ appartient \`a $p^{-1}_{|\hat U_l}(C_l)$, pour $l>k$. On choisit alors $l$ maximal. Par minimalit\'e de $n$ et comme $x$ n'appartient pas \`a $f^{-1}(\cup_{j< k} C_j)$, le point $x':=f^n(x)$ appartient \`a $C_k$. Comme la fibre $F_{x'}^k$ contient $F_{x'}^l$, on a :
\[\ker T_{x'}r_k=\ker T_{x'}r_l\cap \ker T_{x'} (r(z))_{z\in F_{x'}^k\setminus F_{x'}^l}.\]
Par hypoth\`ese de r\'ecurrence, on a :
\[\ker T_{x'}r_k\subset \chi_l(x') \cap \ker T_{x'} (r(z))_{z\in F_{x'}^k\setminus F_{x'}^l}.\]
On va montrer que les \'el\'ements de $F_{x'}^k\setminus F_{x'}^l$ appartiennent \`a $D$. On a alors d'apr\`es $(A_l)$  :
 \[\ker T_{x'}r_k\subset \chi_l(x') \cap \ker T_{x'} (R_1(z))_{z\in F_{x'}^k\setminus F_{x'}^l}\subset \chi_k(x').\]
Et par $f_*$-stabilit\'e de $\chi_k$, on a :
 \[\ker T_{x}r_k\subset (f_*^n\chi_k)(x)\subset \chi_k(x).\]
Ce que l'on souhaitait d\'emontrer.\\
\indent Il suffit donc de montrer que les \'el\'ements de $F_{x'}^k\setminus F_{x'}^l$ appartiennent \`a $D$.  Tout d'abord, le point $x'$ appartient \`a $C_l$. Donc, tous les points de $F_{x'}^k\setminus F_{x'}^l$ n'appartiennent pas \`a 
$ \cup_{j<l}\hat f^{-1}(p^{-1}_{|\hat U_j}(W_j))=\cup_{j<l}\hat f^{-1}(p^{-1}_{|\hat U_j}(C_j))$. Par d\'efinition, ces \'el\'ements n'appartiennent pas n'ont plus \`a $p_{|\hat U_l}^{-1}(C_l)$. Enfin par maximalit\'e de $l$, l'ensemble  $F_{x'}^k\setminus F_{x'}^l$ ne peut pas intersecter $\cup_{j>l}p^{-1}_{|\hat U_j}(C_j)$. Ainsi, l'ensemble $F_{x'}^k\setminus F_{x'}^l$ est inclus dans $\hat f^{-1}(U)= \cup_j\hat f^{-1}(p_{|\hat U_j}^{-1}(C_j))$. Comme $x'$ appartient \`a $C_k$, les \'el\'ements de $F_{x'}^k\setminus F_{x'}^l$ appartiennent bien \`a $D$.
\end{proof}

\indent Cette derni\`ere propri\'et\'e montre en particulier que $\ker(T_xr_k)$ est un espace de dimension $k$, pour tout $x\in C_k$.\\
\indent On montre maintenant par r\'ecurrence d\'ecroissante sur $k$ que la propri\'et\'e $(P_3)$ est v\'erifi\'ee. Soit $k\in\{1,\dots, d\}$. On va commencer par montrer que $\ker (Tr_k)$ est continue sur $C_k$. Par l'hypoth\`ese de r\'ecurrence, seule la continuit\'e en $\tilde K_k:=C_k\cap X_k$ n'est pas \'evidente. Soit $(x_n)_n\in C_k^{\mathbb N}$ une suite qui converge vers $x\in \tilde K_k$. Soit $E$ une valeur d'adh\'erence de $(\ker(T_{x_n}r_k))_n$. L'angle entre les espaces $E$ et $E_k(x)=T_xX_k$ est donc inf\'erieur \`a $\epsilon_k$. Par $f_*$-stabilit\'e de $\ker T r_k$ et $f$-stabilit\'e de $\tilde K$, il existe une valeurs d'adh\'erence $E_m$ de $(\ker(T_{f^{m}(x_n)}r_k))_n$, qui est $\epsilon_k$-proche de $E_k(f^m(x))$ et telle que  $E$ soit \'egal \`a $(T_xf^m)^{-1}(E_m)$. Donc par dilatation normale et la propri\'et\'e 1 de \ref{faitcone}, l'espace $E$ est \'egal \`a $E_k(x)=T_xX_k$. Cela prouve la continuit\'e de $\ker(Tr_k)$, par compacit\'e de la grassmannienne.\\    
\indent On finit maintenant de d\'emontrer la propri\'et\'e $(P_3)$. Pour $x\in X_k$, il existe $n\ge 0$ tel que le point $f^{n}(x)$ appartient \`a l'int\'erieur de $\tilde K_k$ dans $X_k$. Par dilatation normale, il existe un voisinage $V_x$ de $x$ dont l'image par $f^n$ est incluse dans $C_k$ et tel que, pour tout $x'\in V_x$, l'espace $\ker T_{x'}r_k$ est de dimension $k$.\\
\indent Par r\'egularit\'e de $Tf^n$ et r\'egularit\'e de $\ker Tr_{k|C_k}$, quitte \`a r\'eduire $V_x$, la restriction $\ker Tr_{k|V_x}$ est continue.\\
\indent On pose alors $L_k:= int(\cup_{x\in X_k} V_x\cup C_k)$, qui v\'erifie donc la propri\'et\'e $(P_3)$. 
\paragraph{vi Construction de $\Psi$ (preuve de la propri\'et\'e \ref{proprPsi})\\}
La difficult\'e de cette propri\'et\'e \label{Psi}  r\'eside dans le fait que la famille de r\'eels $(t_k)_k$ soit r\'ecursivement assez petite et que les compacts $(C_k)_k$ s'intersectent.\\
\indent On rappelle que, dans la partie i), on a d\'efini les r\'eels $C>\lambda>1$. On fixe une fonction $\phi$ croissante de classe $C^s$ sur $\mathbb R$, s'annulant sur $]-\infty, 1/\lambda]$ et \'egale \`a 1 sur $]1,\infty[$. Pour $t> 0$, on pose $\phi_t:=\phi(\cdot/t)$.\\
\indent Soit 
$\Psi:=\prod_{j=1}^d \phi_j$\quad
avec\quad   $\phi_k:= z\in \Upsilon \mapsto \left\{ \begin{array}{cl}
\phi_{t_k}(\sum_{y\in F^k_{p(z)}} R(y))&\mathrm{si}\; z\in \hat U_k\\
1&\mathrm{sinon}\end{array}\right.$.\\
\indent Remarquons que $\Psi$ est de classe $C^s$ quand $(t_k)_k$ est r\'ecursivement assez petite : pour $k\in\{0,\dots ,d-1\}$, il suffit que $t_k$ soit assez petit pour que l'adh\'erence de $p^{-1}_{|\hat U_k}(W_k)$ intersect\'ee avec la fronti\`ere de $\hat U_k$ soit incluse dans $\hat f^{-1}\big(\cup_{j<k} p^{-1}_{|\hat U_j} (W_j)\big)$, o\`u $\prod_{j=1}^{k-1} \phi_j$ est nulle.\\
\indent On remarque enfin que la fonction $\Psi$ est bien nulle sur $\Upsilon\setminus \hat f^{-1}(U)$ et \'egale \`a 1 sur $U$. \\
\indent On doit donc montrer que, pour $(t_j)_j$ r\'ecursivement assez petite, le noyau $E_k(x):= \ker(dR_1\circ Tp_y)_{y\in F_x^k}$ est uniform\'ement proche de celui de $x\mapsto (dR\circ Tp_y)_{y\in F_x^k}$,  pour $x\in C_k$ et $k\in\{0,\dots ,d-1\}$.\\
\indent Pour toute la suite de cette preuve, les estimations seront uniformes sur $C_k$ ou sur $ p^{-1}_{|\hat U_k}(C_k)$ et seront effectu\'ees pour $(t_k)_k$ r\'ecursivement assez petite.\\
\indent On commence par calculer la diff\'erentielle de $R_1$ :
\[dR_1 = \Psi dR+\frac{(1-\Psi)}{C}dR\circ T\hat f+\left(R- \frac{R\circ \hat f}{C}\right)d\Psi.\]  
 Et, on a pour $z\in \Upsilon$ :
\[d_z\Psi=\sum_{\{i;\; \hat U_i\ni z\}} \left(\prod_{j\not = i} \phi_j\right)(z) \cdot d_z\phi_i
=\sum_{\{i:\; \hat U_i\ni z\}} \underbrace{\left(\prod_{j\not = i} \phi_j\right)(z) \cdot \frac{\phi'}{t_i}}_{=: f_i(z)}\cdot\sum_{y\in F_{p(z)}^i} d(R\circ p_y).\]

\indent On analyse maintenant la diff\'erentielle de $R_1$.\\
\begin{itemize}
\item Les fonctions $\Psi$ et $(1-\Psi)/C$ sont \`a valeurs dans $[0,1]$.
 \item Par dilatation normale, la diff\'erentielle $\frac{dR\circ T\hat f}{\|dR\circ T\hat f\|}$ est proche de $\frac{dR}{\|dR\|}$ sur $p^{-1}_{|\hat U_k}(C_k)$. On a ainsi l'existence d'une fonction  continue $a$ sur $\Upsilon$, born\'ee et sup\'erieure \`a $1$, telle que :
 \[\Psi dR+\frac{(1-\Psi)}{C}dR\circ T\hat f=a\cdot dR+o(1),\quad \mathrm{sur}\; p_{|\hat U_k}^{-1}(C_k).\]
On note que la fonction $a$ est ind\'ependante de $(t_j)_j$.  

\item La fonction $\rho := \left(R-\frac{R\circ \hat f}{C}\right)$ est positive et inf\'erieure \`a $R$, d'apr\`es (i.1). Donc, pour tout $l$,  sur $p_{|\hat U_l}^{-1}(C_l)$, la fonction $\rho$ est \`a valeurs dans $[0,t_l]$.
  \end{itemize}
 
\indent Malheureusement, la norme uniforme de $\rho\cdot d\Psi$ sur $C_k$ n'est ni n\'egligeable ni dans la direction de $dR$, pour $k\in\{1,\dots, d-2\}$. Cependant, la propri\'et\'e \ref{Psi} veut seulement que l'intersection des noyaux de $(dR_1\circ Tp_y)_{y\in F_x^k}$ soit proche de l'intersection des noyaux de $(dR\circ Tp_y)_{y\in F_x^k}$, sur $C_k$.\\     
\indent On remarque que, pour $i<l$, la norme uniforme de la fonction $f_i$  est petite devant $1/t_l$.\\ 
Ainsi, pour $x\in C_l$ et $z\in F_x^l$, on a : 
\[\rho(z)\cdot \sum_{\{i<l:\; \hat U_i\ni z\}} f_i(z)\sum_{y\in F_x^i} d(R\circ p_y)=o(1).\]
Et, pour $z\in \hat U_i\setminus p_{|\hat U_i}^{-1}(C_i)$, on a :
\[f_i(z)=0\]
\indent On conclut que, pour $x\in C_k$ et $z\in F_x^k$, si $l\ge k$ est maximal tel que $z$ soit dans $p_{|\hat U_l}^{-1}(C_l)$, on a :  
\[d_zR_1= a(z)\cdot d_zR + \rho(z)\cdot  f_{l}(z) \sum_{y\in F_x^{l}} d(R\circ p_y) +o(1).\]

\indent Aussi, pour $k\in\{1,\dots, d\}$ et $x\in C_k$, si $x$ appartient exactement \`a $(C_{i_j})_{j=1}^l$, pour $(i_j)_j\in\{k,\dots, d-1\}^l$ d\'ecroissant (et ainsi $i_l=k$), on a pour $z\in F_x^{i_j}\setminus F_x^{i_{j-1}}$ (avec $F_x^0:=\emptyset$) :
\[d_zR_1= a(z)\cdot d_zR + \rho(z)\cdot  f_{i_j}(z) \sum_{y\in F_x^{i_j}} d(R\circ p_y) +o(1).\]
On munit $F_x^k$ d'un ordre compatible avec l'indexation $(i_j)_j$, selon lequel on effectue un produit ext\'erieur :
\[\bigwedge_{z\in F_x^k}d(R_1\circ p_z)=\bigwedge_{j=1}^l\bigwedge_{z\in F_x^{i_j}\setminus F_x^{i_{j-1}}}\left( a(z)d(R\circ p_z)+ \rho(z)\cdot f_{i_j}(z) \sum_{y\in F_x^{i_j}} d(R\circ p_y) +o(1)\right).\]

Tous les scalaires \'etant positifs, ce produit est \'egal \`a :
\[ \prod_{j=1}^l\left(\prod_{z\in F_x^{i_j}\setminus F_x^{i_{j-1}}} a(z)+\sum_{z\in F_x^{i_j}\setminus F_x^{i_{j-1}}} \rho(z)\cdot  f_{i_j}(z)\prod_{y\in F_x^{i_j}\setminus (F_x^{i_{j-1}}\cup \{z\})} a(y)\right)\bigwedge_{z\in F_x^k}\Big(d(R\circ p_z)+o(1)\Big).\]

Cela implique le noyau de $(d(R_1\circ p_y))_{y\in F_x^k}$, est de dimension $k$ et uniform\'ement proche de $\ker(d(R\circ p_y))_{y\in F_x^k}$ sur $C_k$, pour une famille $(t_j)_j$ r\'ecursivement assez petite.
\begin{flushright}
$\square$
\end{flushright}
\subsection{Preuve de la propri\'et\'e $\ref{4.2.2}$}\begin{proof}[] Pour chaque $k\ge 0$, comme $N$ est suppos\'ee envoy\'ee par $f$ dans $N'$ relativement compacte, par dilatation normale,  il existe un compact $K$ de $X_k$ tel que  
\[\cup_{n\ge 0} f_{|N}^{-n}(K)=X_k.\]
Comme $p_{|A_k}$ est un rev\^etement \`a $d-k$ feuillets de $X_k$ et comme $p$ est un diff\'eomorphisme local, il existe un voisinage ouvert $\hat V_k$ de $p_{|A_k}^{-1}(K)$ tel que $p_{|\hat V_k}$ et $p_{|\hat V_k \cup \hat f^{-1}(\hat V_k)}$ soient des rev\^etements \`a $d-k$ feuillets de $V_k:= p(\hat V_k)$ et $V_k\cup f_{|N}^{-1}(V_k)$ respectivement.
\[ \mathrm{ Soient}\;  \hat U_k:= \cup_{n\ge 0} \hat f^{-n}(\hat V_k)\; \mathrm{et} \; U_k:=p(\hat U_k).\]
\indent D'apr\`es l'expression de $\hat U_k$, la pr\'eimage $\hat f^{-1}(\hat U_k)$ est incluse dans $\hat U_k$.  On va montrer que $f^{-1}(U_k)$ est inclus dans $U_k$.\\ %On remarque tout d'abord que $\hat f$ pr\'eserve $A_k$. 
\indent On consid\`ere l'application  $I\; :\; x\in \hat V'\mapsto (p(x),\pi(x))\in V'\times \partial N$, avec $\pi\; :\; \partial N\times \mathbb R^+\rightarrow \partial N$ la projection canonique. L'application $I$ est une bijection qui permet de faire l'identification g\'eom\'etrique suivante.  Un point de $\hat V'$ la donn\'ee d'un point $z$ de $N$ avec une face de l'intersection de $N$ avec un voisinage ouvert de $z$. Appelons cette face, \emph{une face locale de $N$ proche de $x$}. Les faces locales ne doivent pas \^etre confondues avec les faces (globales) de $N$.
Par exemple, quand la vari\'et\'e est la "goutte" (voire la partie \ref{rappel}), un point proche de l'arr\^ete de dimension $0$ est proche de deux faces locales, alors que cette vari\'et\'e est munie d'une unique face. \\
% Plus formellement, cette face locale correspond \`a l'ensemble des petits voisinages ouverts  de la projection canonique d'un certain $z'\in p^{-1}(\{z\})\in \partial N\times \mathbb R^+$ dans $\partial N$.\\  
\indent Par dilatation normale, en choisissant $\hat V_k$ petit, on peut choisir $\hat U_k$ inclus dans un petit voisinage de l'adh\'erence de $A_k$.\\ 
\indent L'endomorphisme $Tf$ induit une application entre les faces locales proches de $f(x)\in U_k$ vers les faces locales proches de $x$ car, par dilation normale, la pr\'eimage par $Tf$ d'un hyperplan proche de $TX_{d-1}$ est un hyperplan proche de $TX_{d-1}$. La dilatation normale entraine m\^eme que cette application entre faces locales est bijective. On remarque que l'inverse de cette application est $\hat f$. Ainsi, pour $x\in f^{-1}(U_k)$, un point $y\in F_{f(x)}^k$ est la donn\'ee de $f(x)$ munie d'une face locale proche de $f(x)$. Cette face poss\`ede une pr\'eimage point\'ee en $x$. Cela entra�ne que $y$ a une pr\'eimage par $\hat f$. Ainsi $x$ appartient \`a $p(\hat f^{-1}(\hat U_k))$ qui est inclus dans $U_k=p(\hat U_k)$.

\indent Cela entra�ne aussi que, pour $x\in f^{-1}_{|N}( U_k)$, la fibre $F_x^k$ est incluse dans $\hat f^{-1}(\hat U_k)$ et que $\hat f_{|F_x^k}$ est une bijection de $F_x^k$ sur $F_{f(x)}^k$.

On va  maintenant montrer les \'egalit\'es suivantes :
\begin{equation}\label{hatuk}U_k=\cup_{n\ge 0}f^{-n}(V_k)\quad \mathrm{et} \quad p(\hat f^{-n}(\hat V_k'))= f^{-n}(V_k'),\end{equation}\[ \mathrm{avec} \quad \hat V_k':=\hat f^{-1}(\hat V_k)\setminus  \hat V_k\quad\mathrm{et}\quad
 V_k':= p(\hat V_k')=f_{|N}^{-1}(V_k)\setminus V_k.\]
\indent  Soient $x\in U_k\setminus V_k$ et $y\in F_x^k$, comme $\hat U_k$ est \'egal \`a $\hat V_k\cup \bigcup_{n\ge 0} \hat f^{-n}(\hat V_k')$, il existe alors un unique $n\ge 0$ tel que $y$ appartient \`a $\hat f^{-n}(\hat V_k')$. Donc $\hat f^n(y)$ appartient \`a $\hat V_k'$ et, par commutativit\'e du diagramme, $f^n(x)$ appartient \`a $V_k'$. Cela implique que $x$ appartient \`a $f^{-n}(V_k')$. Donc $p(\hat f^{-n}(\hat V_k'))$ est inclus dans  $f^{-n}(V_k')$ et $U_k$ est inclus dans $\cup_{\ge 0} f^{-n}(V_k)$. Comme $U_k$ contient $V_k$ et est $f^{-1}$ stable, on obtient les  \'egalit\'es de (\ref{hatuk}).

On peut maintenant d\'emontrer que $p_{|\hat U_k}$ est un rev\^etement \`a $(d-k)$ feuillets de $U_k$. On a d\'efini $\hat V_k$ de fa\c con \`a ce que la restriction de $p$ \`a $\hat V_k\cup \hat f^{-1}(V_k)=\hat V_k\cup \hat V_k'$ soit un rev\^etement \`a $(d-k)$ feuillets de $V_k\cup f^{-1}(V_k)=V_k\cup V_k'$. Par la deuxi\`eme \'egalit\'e de (\ref{hatuk}), tous les \'el\'ements de $\hat U_k\setminus (\hat V_k\cup \hat V_k')$ sont envoy\'es par $p$ \`a l'ext\'erieur de $V_k\cup V_k'$. Donc $F_x^k$ est de cardinalit\'e $(d-k)$ pour $x\in V_k\cup V_k'$.\\
\indent Pour $x\in f^{-n}(V_k')$, par la deuxi\`eme \'egalit\'e de (\ref{hatuk}), la fibre $F_x^k$ est incluse dans $\hat f^{-n}(\hat V_k)$. Par la bijectivit\'e de $\hat f^{n}_{|F_x^k}$, puis la commutativit\'e du diagramme, et enfin la cardinalit\'e de $F^k_y$, pour $y=f^n(x)\in V_k'$, la fibre $F_x^k$ est bien de cardinalit\'e $(d-k)$.\\
\indent D'apr\`es la premi\`ere \'egalit\'e de (\ref{hatuk}), cette discussion implique que $p_{|\hat U_k}$ est un rev\^etement de $U_k$ \`a $(d-k)$ feuillets.

\indent Pour $j\le k$ et $x\in X_k$ qui tend vers $y\in X_j$, $F_{x}^k$ tend \`a \^etre inclus dans $F^j_{y}$. Quitte \`a restreindre $V_k$ et $V_j$, on peut donc supposer que, pour $x\in U_k\cap U_j$, la fibre $F_{x}^k$ est incluse dans $F_{x}^j$.\end{proof}%\\} 

\section{Questions et remarques sur le th\'eor\`eme de persistance des vari\'et\'es \`a coins}
En r\'ealit\'e, la version g\'en\'erale du th\'eor\`eme \ref{thfonda} (voir le r\'esultat principal de \cite{PB1}) implique que nous avons montr\'e la $C^s$-persistance des vari\'et\'es \`a coins en un sens plus fort. En effet, toute perturbation $f'$ de $f$ pr\'eserve une stratification $\Sigma'$ qui est l'image de $\Sigma_{|N'}$ par un plongement $p$ contr\^ol\'e :
\begin{itemize} 
\item $p$ est un hom\'eomorphisme sur son image $C^0$-proche de l'inclusion canonique de $N'$ dans $M$,
\item $p$ envoie les strates de $\Sigma_{|N'}$ sur les strates de $\Sigma'$ qui sont pr\'eserv\'ees par $f'$.
\item la restriction de $p$ \`a chaque intersection de $N'$ avec chaque lamination\footnote{$(L_k,\mathcal L_k)$ est d\'efinie durant la preuve du th\'eor\`eme \ref{varcoin}.}  $(L_k, \mathcal L_k)$ est un plongement de lamination de classe $C^s$, proche de l'inclusion canonique. Cela signifie en plus que :
\begin{itemize}\item  $p$ est un plongement de classe $C^s$ le long des plaques de $\mathcal L_k$ contenues dans $N'$, \item  ses $s$ premi\`eres diff\'erentielles varient continument transversalement au plaques,
\item ses diff\'erentielles sont proches de celle de l'identit\'e sur tout compact de $L_k$,\end{itemize} \end{itemize}
Ainsi, $\Sigma'$ est munie d'une structure de treillis de laminations $\mathcal T'$ de classe $C^s$. La version g\'en\'erale du th\'eor\`eme \ref{thfonda} implique enfin que la structure de treillis $\mathcal T'$ v\'erifie les hypoth\`eses $(i)$ et $(ii)$ du th\'eor\`eme \ref{thfonda}.

Chacune des laminations $(L_k,\mathcal L_k)_k$ \'etant des fibrations, il semble facile de constater que la stratification $\Sigma_{|N'}$ est persistante en tant que stratification $c$-r\'eguli\`ere au sens de Bekka (voir \cite{Bek}).         
\paragraph{Questions}  \begin{itemize}
\item R\'eciproquement, on peut se demander si toute stratification $c$-r\'eguli\`ere, pr\'eserv\'ee et normalement dilat\'ee par la dynamique est persistante en tant que stratification $c$-r\'eguli\`ere.% Pour ce r\'esultat, il manque essentiellement la preuve de la conjecture \ref{trot} de Trotman pour les stratifications $c$-r\'eguli\`eres. 
\item Les orbifolds g\'en\'eralisent les structures de vari\'et\'es \`a coins et d\'efinissent aussi canoniquement des stratifications. Consid\'erons un orbifold plong\'e dans une vari\'et\'e, dont la stratification canonique est normalement dilat\'ee. cette stratification est-elle persistante?
\item De quelle fa\c con les sous-vari\'et\'es \`a bord ou \`a coins, analytiques r\'eelles ou complexes, persistent? 
\item Ma\~ne a montr\'e que les sous-vari\'et\'es compactes $C^1$-persistantes et uniform\'ement localement maximales \cite{Manethese} sont les sous vari\'et\'es normalement hyperboliques. Hirsh-Pugh-Shub \cite{HPS} ont montr\'e que toute sous-vari\'et\'e normalement hyperbolique, par un diff\'eomorphisme $f$,  est l'intersection transverse de deux sous vari\'et\'es (une fois) normalement dilat\'ees et dont les adh\'erences compactes sont envoy\'ees dans elles m\^eme par respectivement $f$ et $f^{-1}$. Est ce que toute vari\'et\'e \`a coins compacte, $C^1$-persistante (en tant que stratification) et uniform\'ement localement maximale est aussi l'intersection transverse de deux sous-vari\'et\'es \`a coins v\'erifiant notre th\'eor\`eme \ref{varcoin} pour respectivement $f$ et $f^{-1}$ ? \end{itemize}
 R\'eciproquement, la persistance des vari\'et\'es \`a coins de "fa\c con contr\^ol\'ee" implique que toute vari\'et\'e \`a coins compact, laiss\'ee invariante par un diff\'eomorphisme $f$, qui est l'intersection transverse de deux vari\'et\'es \`a coins v\'erifiant notre th\'eor\`eme \ref{varcoin} pour respectivement $f$ et $f^{-1}$,  est alors persistante en tant que stratification $a$-r\'eguli\`ere.    
\bibliographystyle{alpha}
\nocite{*}
\bibliography{references}

\def\polhk#1{\setbox0=\hbox{#1}{\ooalign{\hidewidth
  \lower1.5ex\hbox{`}\hidewidth\crcr\unhbox0}}} \def\cprime{$'$}
\begin{thebibliography}{HPS77}

\bibitem[Bek91]{Bek}
K.~Bekka.
\newblock C-r\'egularit\'e et trivialit\'e topologique.
\newblock In {\em Singularity theory and its applications, Part I (Coventry,
  1988/1989)}, volume 1462 of {\em Lecture Notes in Math.}, pages 42--62.
  Springer, Berlin, 1991.

\bibitem[Ber]{PB1}
Pierre Berger.
\newblock Persistence of stratification of normally expanded laminations.
\newblock {\em arXiv:math.DS}.

\bibitem[Cer61]{C}
Jean Cerf.
\newblock Topologie de certains espaces de plongements.
\newblock {\em Bull. Soc. Math. France}, 89:227--380, 1961.

\bibitem[Dou62]{D}
Adrien Douady.
\newblock Vari\'et\'es \`a bord anguleux et voisinages tubulaires.
\newblock In {\em S\'eminaire Henri Cartan, 1961/62, Exp. 1}, page~11.
  Secr\'etariat math\'ematique, Paris, 1961/1962.

\bibitem[Hir76]{H}
Morris~W. Hirsch.
\newblock {\em Differential topology}.
\newblock Springer-Verlag, New York, 1976.
\newblock Graduate Texts in Mathematics, No. 33.

\bibitem[HPS77]{HPS}
M.~W. Hirsch, C.~C. Pugh, and M.~Shub.
\newblock {\em Invariant manifolds}.
\newblock Springer-Verlag, Berlin, 1977.
\newblock Lecture Notes in Mathematics, Vol. 583.

\bibitem[Ma{\~n}78]{Manethese}
Ricardo Ma{\~n}{\'e}.
\newblock Persistent manifolds are normally hyperbolic.
\newblock {\em Trans. Amer. Math. Soc.}, 246:261--283, 1978.

\bibitem[Mat73]{Ma}
John~N. Mather.
\newblock Stratifications and mappings.
\newblock In {\em Dynamical systems (Proc. Sympos., Univ. Bahia, Salvador,
  1971)}, pages 195--232. Academic Press, New York, 1973.

\bibitem[Shu69]{Shubthese}
Michael Shub.
\newblock Endomorphisms of compact differentiable manifolds.
\newblock {\em Amer. J. Math.}, 91:175--199, 1969.

\bibitem[Tho64]{Th}
R.~Thom.
\newblock Local topological properties of differentiable mappings.
\newblock In {\em Differential Analysis, Bombay Colloq.}, pages 191--202.
  Oxford Univ. Press, London, 1964.

\bibitem[Whi65]{W1}
Hassler Whitney.
\newblock Local properties of analytic varieties.
\newblock In {\em Differential and Combinatorial Topology (A Symposium in Honor
  of Marston Morse)}, pages 205--244. Princeton Univ. Press, Princeton, N. J.,
  1965.

\end{thebibliography}
\noindent------------------------\\
\noindent Pierre Berger\\
Institute for Mathematical Sciences, Stony Brook University, Stony Brook NY 11794-3660, USA (pierre.berger((at))normalesup.org)\\
 
\end{document}